\documentclass[11pt]{amsart}

\usepackage{bm,url,amssymb,amscd}
\usepackage{graphics,epsfig,color}
\usepackage[curve]{xypic}

\def\e{{\rm  e}}
\def\A{{\mathcal A}}
\def\C{{\mathbb C}}

\def\deg{{\rm deg}}

\renewcommand{\em}{\bf}
\renewcommand{\bar}{\overline}
\renewcommand{\Re}{{\rm Re}}
\renewcommand{\Im}{{\rm Im}}

\newtheorem{definition}{\bf Definition}
\newtheorem*{theorem*}{\bf Theorem}
\newtheorem*{TOP_CHAR_Mtheorem}{\bf Topological Characterization of the Mandelbrot Set\index{Topological Characterization of the Mandelbrot Set}}
\newtheorem*{TOP_CHAR_Mtheorem_Prime}{\bf Topological Characterization of the Mandelbrot Set'}

\newtheorem*{TRIANGLE_INEQUALITY}{\bf Triangle Inequality}
\newtheorem*{fatoujulialemma}{\bf Fatou-Julia Lemma\index{Fatou-Julia Lemma}}
\newtheorem*{koenigstheorem}{\bf K\oe nig's Theorem\index{K\oe nig's Theorem}}
\newtheorem*{bottcherstheorem}{\bf B\"ottcher's Theorem\index{B\"ottcher's Theorem}}
\newtheorem*{FTA}{\bf Fundamental Theorem of Algebra\index{Fundamental Theorem of Algebra}}
\newtheorem*{HBtheorem}{\bf Heine-Borel Theorem\index{Heine-Borel Theorem}}
\newtheorem*{IFtheorem}{\bf Inverse Function Theorem\index{Inverse Function Theorem}}
\newtheorem*{SCIFtheorem}{\bf Simply-Connected Inverse Function Theorem\index{Simply-Connected Inverse Function Theorem}}
\newtheorem*{CRtheorem}{\bf Cauchy-Riemann Equations\index{Cauchy-Riemann Equations}}
\newtheorem*{Ctheorem}{\bf Cauchy's Theorem\index{Cauchy's Theorem}}
\newtheorem*{Mtheorem}{\bf Morera's Theorem}
\newtheorem*{CIFtheorem}{\bf Cauchy Integral Formula\index{Cauchy Integral Formula}}
\newtheorem*{CIF2theorem}{\bf Cauchy Integral Formula For Higher Derivative\index{Cauchy Integral Formula For Higher Derivatives}}
\newtheorem*{CEtheorem}{\bf Cauchy Estimates\index{Cauchy Estimates}}
\newtheorem*{MMtheorem}{\bf Maximum Modulus Principle\index{Maximum Modulus Principle}}
\newtheorem*{PStheorem}{\bf Existence of Power Series}
\newtheorem*{ULAtheorem}{\bf Uniform Limits Theorem\index{Uniform Limits Theorem}}
\newtheorem*{ULHtheorem}{\bf Uniform Limits of Harmonic Functions\index{Uniform Limits of Harmonic Functions}}
\newtheorem*{Slemma}{\bf Schwarz Lemma\index{Schwarz Lemma}}
\newtheorem*{APPlemma}{\bf Attracting Periodic Orbit Lemma\index{Attracting Periodic Orbit Lemma}}
\newtheorem*{MAXtheorem}{\bf Maximum Principle\index{Maximum Principle}}
\newtheorem*{AP}{\bf Argument Principle\index{Argument Principle}}
\newtheorem{lemma}{\bf Lemma}
\newtheorem*{corollary}{\bf Corollary}

\newtheorem*{conjecture*}{\bf Conjecture}
\newtheorem*{DENSITYconjecture}{\bf Density of Hyperbolicity Conjecture\index{Density of Hyperbolicity Conjecture}}
\newtheorem*{MLCconjecture}{\bf MLC Conjecture\index{MLC Conjecture}}

\theoremstyle{remark}
\newtheorem{exercise}{\bf Exercise}
\newtheorem{extendedexercise}{\bf Extended Exercise}

\newtheorem*{question*}{\bf Question}
\newtheorem*{remark}{\bf Remark}
\newtheorem*{EandL}{\bf Exponential and Logarithm}
\newtheorem{example}{\bf Example}

\begin{document}

\author{Roland K.W. Roeder}

\begin{abstract}
We introduce the exciting field of complex dynamics at an undergraduate level
while reviewing, reinforcing, and extending the ideas learned in an typical
first course on complex analysis.  Julia sets and the famous Mandelbrot set
will be introduced and interesting properties of their boundaries will be
described.  We will conclude with a discussion of problems at the boundary
between complex dynamics and other areas, including a nice application of the
material we have learned to a problem in astrophysics.
\end{abstract}

\title{Around the boundary of complex dynamics. }

\maketitle


\section*{Preface}

These notes were written for the 2015 Thematic Program on 
Boundaries and Dynamics held at Notre Dame University.  They are intended for an advanced
undergraduate student who is majoring in mathematics.
In an ideal world, a student reading
these notes will have already taken undergraduate level courses in complex
variables, real analysis, and topology.  As the world is far from ideal, we
will also review the needed material.

There are many fantastic places to learn complex dynamics, including
the books by Beardon~\cite{BEARDON}, Carleson-Gamelin \cite{CG}, Devaney \cite{DEV1,DEV2}, Milnor
\cite{MILNOR}, and Steinmetz \cite{STEINMETZ}, as well as the Orsay Notes\index{Orsay Notes} \cite{ORSAY} by Douady\index{Douady} and Hubbard\index{Hubbard}, the surveys by Blanchard \cite{BLANCHARD} and Lyubich \index{Lyubich} \cite{LYUBICH_SURVEY1,LYUBICH_SURVEY2}, and the invitation to transcendental dynamics by Shen and Rempe-Gillen \cite{SR}.  The books by Devaney and the article by Shen and Rempe-Gillen are especially accessible to undergraduates.
We will take a complementary approach, following a somewhat
different path through some of the same material as presented in these sources.
We will also present modern connections at the boundary between complex
dynamics and other areas.

None of the results presented here are new.  In fact, I learned most of them
from the aforementioned textbooks and from courses and informal discussions
with John Hubbard and Mikhail Lyubich.

Our approach is both informal and naive.  We make no effort to provide a
comprehensive or historically complete introduction to the subject.  Many important results
will be omitted.  Rather, we will simply have fun doing mathematics.

\vspace{0.1in}
\begin{center}
{\it Dedicated to Emile and Eli.}
\end{center}

\vspace{0.1in}
\noindent
{\bf Acknowledgments}
I am grateful to Notre Dame University for their hospitality during the
thematic program on boundaries and dynamics.  Ivan Chio, Youkow Homma, Lyndon
Ji, Scott Kaschner, Dmitry Khavinson, Seung-Yeop Lee,  Rodrigo P\'erez, and
Mitsuhiro Shishikura provided many helpful comments.  All of the computer-drawn images of basins of attraction, filled Julia sets, and the Mandelbrot set were created using the Fractalstream software \cite{FS} that was written by Matthew Noonan.  This work was partially
supported by NSF grant DMS-1348589.\index{computer software}

\section*{Lecture 1: ``Warm up''}
\setcounter{section}{1}
\setcounter{subsection}{0}

Let us start at the very beginning:

\subsection{Complex Numbers}
Recall that a complex number has the form $z = x + i y$, where $x,y \in
\mathbb{R}$ and $i$ satisfies $i^2 = -1$.  One adds, subtracts, multiplies, and divides complex
numbers using the following rules:
\begin{align*}
(a+bi) \pm (c+di) &= (a\pm c) + (b \pm d)i, \\
(a+bi) (c+di) &= ac + adi + bci + bdi^2 = (ac-bd) + (ad+bc) i, \quad \mbox{and}\\
\frac{a+bi}{c+di} &= \frac{a+bi}{c+di}\frac{c-di}{c-di} = \frac{(ac+bd) + (bc-ad)i}{c^2+d^2}.
\end{align*}
The set of complex numbers forms a field  $\C$ under the operations of addition and multiplication.

\begin{figure}[h!]
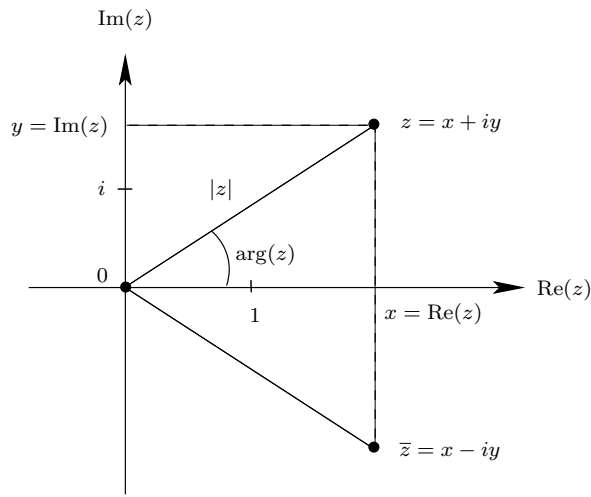
\caption{The complex plane.  \label{FIG_COMPLEX_PLANE}}
\end{figure}

The {\em real part} of $z = x+iy$ is $\Re(z) = x$ and the {\em imaginary part}
of $z = x+iy$ is $\Im(z) = y$.  One typically depicts a complex number in the
{\em complex plane} using the horizontal axis to measure the real part and the
vertical axis to measure the imaginary part; See Figure \ref{FIG_COMPLEX_PLANE}.
One can also take the real or imaginary part of more complicated expressions.  For example,
$\Re(z^2) = x^2-y^2$ and $\Im(z^2) = 2xy$.

The {\em complex conjugate} of $z=x+iy$ is $\overline{z} = x-iy$ and the {\em
modulus} of $z$ is $|z| = \sqrt{x^2+y^2} = \sqrt{z \bar{z}}$.  In the complex
plane, $\bar{z}$ is obtained by reflecting $z$ across the real axis and $|z|$
is the distance from $z$ to the origin $0 = 0+0i$.  The {\em  argument}
of $z \neq 0$ is the angle counterclockwise from the positive real axis to $z$.

A helpful tool is the:
\begin{TRIANGLE_INEQUALITY}
For every $z,w \in \C$ we have
\begin{align*}
|z| - |w| \leq |z+w| \leq |z|+|w|.
\end{align*}
\end{TRIANGLE_INEQUALITY}

A {\em complex polynomial\index{complex polynomial}} $p(z)$ of degree $d$ is an expression of the form
\begin{align*}
p(z) = a_d z^d + a_{d-1} z^{d-1} + \cdots + a_1 z + a_0
\end{align*}
where $a_d, \ldots, a_0$ are some given complex numbers with $a_d \neq 0$.  
Historically, complex numbers were introduced so that the following theorem holds:

\begin{FTA}
A polynomial $p(z)$ of degree $d$  has $d$ complex zeros $z_1, \ldots, z_d$, counted with multiplicity.
\end{FTA}

In other words, a complex polynomial\index{complex polynomial} $p(z)$ can be factored over the complex numbers as
\begin{align}\label{EQN_FACTORIZED_POLY}
p(z) = c (z-z_1)(z-z_2)\cdots(z-z_d),
\end{align}
where $c \neq 0$ and some of the roots $z_j$ may be repeated.  (The number of times $z_j$ is repeated in (\ref{EQN_FACTORIZED_POLY}) is the {\em multiplicity} of $z_j$ as a root of $p$.)

\begin{figure}[h!]
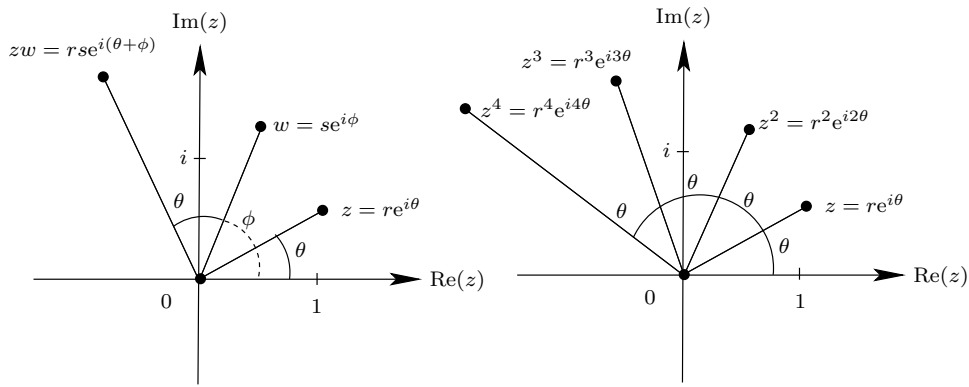
\caption{Multiplication and taking powers in polar form.  \label{FIG_MULT_POWER}}
\end{figure}

Multiplying and dividing complex numbers is often simpler in {\em polar form}. Euler's Formula states 
\begin{align*}
{\rm e}^{i\theta} = \cos \theta + i \sin \theta \qquad \mbox{ for any $\theta \in \mathbb{R}$}.
\end{align*}
We can therefore represent any complex number $z = x+iy$ by $z = r {\rm e}^{i\theta}$ where $r= |z|$
and $\theta = \arg(z)$.  Suppose $z = r \e^{i\theta}$ and $w = s \e^{i\phi}$ and $n \in \mathbb{N}$. The simple formulae
\begin{align}\label{EQN_MULT_DIVIDE_POLAR}
zw = rs \e^{i(\theta+\phi)}, \qquad z^n = r^n  \e^{i n\theta}, \quad \mbox{and} \quad \frac{z}{w} = \frac{r}{s}  \e^{i(\theta-\phi)}.
\end{align}
follow from the rules of exponentiation.
Multiplication and taking powers of complex numbers in polar form are depicted geometrically in Figure \ref{FIG_MULT_POWER}.

\subsection{Iterating Linear Maps} \label{SUBSEC_LINEAR_MAPS}

A {\em linear map} $L: \C \rightarrow \C$ is a mapping of the form $L(z) = az$, where $a \in \C \setminus \{0\}$.
Suppose we take some {\em initial condition\index{initial condition}} $z_0 \in \C$ and repeatedly apply $L$:
\begin{align}\label{EQN_ITERATES_OF_L}
\xymatrix{
z_0 \ar[r] & L(z_0) \ar[r] & L(L(z_0)) \ar[r] & L(L(L(z_0))) \ar[r] & \cdots .
}
\end{align}
For any natural number $n \geq 1$ let
$L^{\circ n}: \C \rightarrow \C$ denote the composition of $L$ with itself $n$ times.
We will often also use the notation
\begin{align*}
z_n := L^{\circ n}(z_0).
\end{align*}
The sequence $\{z_n\}_{n=0}^\infty \equiv \{L^{\circ n}(z_0)\}_{n=0}^\infty$ is called the { sequence of iterates} of $z_0$ under $L$.  It is also called the {\em orbit\index{orbit}} of $z_0$ under $L$.

\begin{remark} The notion of {\em linear} used above is from your course on linear algebra: a linear 
map must satisfy $L(z+w) = L(z) + L(w)$ for all $z,w \in \C$ and $L(cz) = c L(z)$ for all $z,c \in \C$.  For
this reason, mappings of the form $z \mapsto a z + b$ are {\em not considered} linear.  Instead, they are called {\em affine}.
(See Exercise \ref{EX_AFFINE}.)
\end{remark}

The number $a$ is called a {\em parameter\index{parameter}} of the system.  We think of it as
describing the overall state of the system (think, for example, temperature or barometric
pressure) that is fixed for all iterates~$n$.  One can change the parameter\index{parameter} to
see how it affects the behavior of sequences of iterates (for example, if the temperature is higher, does the orbit\index{orbit} move farther in each step?).

\begin{figure}[h!]
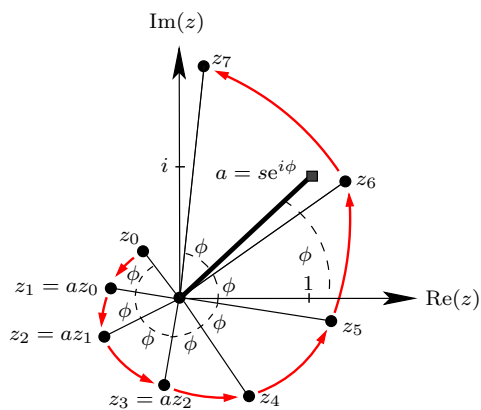
\caption{Iterating the linear map $L(z) = az$.  Above: $|a| < 1$ implies orbits\index{orbit}
spiral into $0$.  Below: $|a| > 0$ implies \index{orbit}s spiral away from $0$.  Not Shown:
$|a| = 1$ implies orbits\index{orbit} rotate around $0$ at constant modulus.
\label{FIG_ITERATING_LINEAR_MAP}}
\end{figure}

Our rules for products and powers in polar form (\ref{EQN_MULT_DIVIDE_POLAR})
allow us to understand the sequence of iterates (\ref{EQN_ITERATES_OF_L}).
Suppose $z_0 = r \e^{i\theta}$ and $a = s \e^{i\theta}$ with $r,s > 0$.  Then,
the behavior of the iterates depends on $s = |a|$, as shown in Figure
\ref{FIG_ITERATING_LINEAR_MAP}.

\begin{remark}
For a linear map $L(z) = az$ with $|a| \neq 1$ the orbits\index{orbit} $\{z_n\}$ and
$\{w_n\}$ for any two non-zero initial conditions\index{initial condition} $z_0$ and $w_0$ have the same
dynamical behavior.  If $|a| < 1$ then $$\lim_{n\rightarrow \infty} z_n = 0 =
\lim_{n\rightarrow \infty} w_n$$ and if $|a| > 1$ then $$\lim_{n\rightarrow \infty}
z_n = \infty = \lim_{n\rightarrow \infty} w_n.$$  

This is atypical for dynamical systems---the long term behavior of the
orbit\index{orbit} usually depends greatly on the initial condition\index{initial condition}.   For example, we will
soon see that when iterating the quadratic mapping $p(z) = z^2 + \frac{i}{4}$ there are many initial conditions\index{initial condition} whose orbits\index{orbit} remain bounded and many whose
orbit\index{orbit} escapes to $\infty$.
There will also be many initial conditions\index{initial condition} whose orbits\index{orbit} have completely different behavior!
Linear maps are just too simple to have interesting dynamical properties.
\end{remark}

\begin{exercise}\label{EX_AFFINE} 
An affine mapping $A: \C \rightarrow \C$ is a mapping given by $A(z) = az + b$, where $a,b \in \C$ and $a \neq 0$.
Show that iteration of affine mappings produces no dynamical behavior that was not seen when iterating linear mappings.
\end{exercise}

\subsection{Iterating quadratic polynomials.} \label{SUBSEC_QUADRATIC_MAPS}

Matters become far more interesting if one iterates quadratic mappings $p_c: \C
\rightarrow \C$ given by $p_c(z) = z^2 + c$.  Here, $c$ is a parameter\index{parameter}, which
we sometimes include in the notation by means of a subscript, writing $p_c(z)$, and sometimes omit, writing simply~$p(z)$.

\begin{remark} Like in Exercise \ref{EX_AFFINE}, one can show that quadratic mappings of the form $p_c(z) = z^2 + c$
actually capture all of the types of dynamical behavior that can arise when iterating a more general quadratic mapping $q(z) = az^2+bz+c$.
\end{remark}

\begin{figure}[h!]
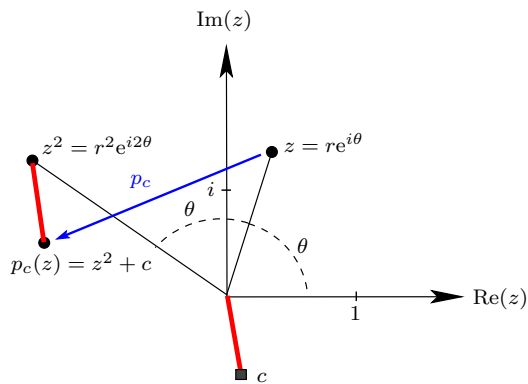
\caption{Geometric interpretation of applying $p_c(z) = z^2+c$.
\label{FIG_APPLYING_QUADRATIC_MAP}}
\end{figure}

Applying the mapping $p_c$ can be understood geometrically in two steps: one first squares
$z$ using the geometric interpretation provided in polar coordinates
(\ref{EQN_MULT_DIVIDE_POLAR}).  One then translates (shifts) the result by $c$.
This two-step process is illustrated in 
Figure~\ref{FIG_APPLYING_QUADRATIC_MAP}.

\begin{remark} Solving the exercises in this subsection may require some of the basic complex analysis from the following subsection.  They
are presented here for better flow of the material.
\end{remark}

\begin{example}{\bf Exploring the dynamics of $\bm p_c: \C \rightarrow \C$ for $\bm c = \frac{i}{4}$.}\label{EXAMPLE_C_I_OVER_4}
In Figure \ref{FIG_ZSQUARED_PLUS_I_OVER_4} we show the first few iterates under $p(z) = z^2+\frac{i}{4}$ of
two different orbits\index{orbit}: $\{z_n\}$ of initial condition\index{initial condition} $z_0 = i$ and $\{w_n\}$ of initial condition\index{initial condition} $w_0~=~1.1i$.
Note that orbit $\{z_n\}$ seems to converge
to a point \hbox{$z \approx -0.05+.228i$} while orbit\index{orbit} $\{w_n\}$ seems to escape to $\infty$.

\begin{figure}[h!]
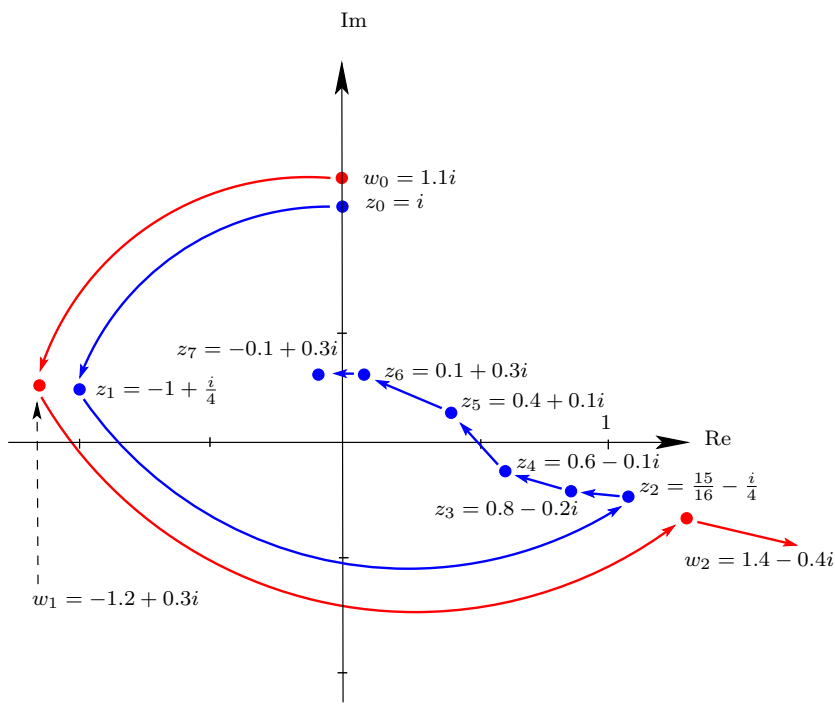\caption{Orbits $\{z_n\}$ for initial condition\index{initial condition} $z_0 = i$ and $\{w_n\}$ for $w_0 = 1.1i$ under $p(z) = z^2+\frac{i}{4}$.
\label{FIG_ZSQUARED_PLUS_I_OVER_4}}
\end{figure}

\begin{exercise}\label{EX_FIXED_POINTS}
Use the quadratic formula to prove that there exists $z_\bullet \in \mathbb{C}$ that is close to $-0.05+.228i$ and satisfies
\begin{align*}
p(z_\bullet) = z_\bullet.
\end{align*}
Such a point is called a {\em fixed point\index{fixed point}} for
$p(z)$ because if you use $z_\bullet$ as initial condition\index{initial condition} the orbit\index{orbit} is a constant sequence $\{z_\bullet,z_\bullet,z_\bullet,\ldots \}$.  

Show that there is a second fixed point\index{fixed point} $z_*$ for
$p(z)$ with $z_* \approx 1.05-.228i$. 

 Compute
$|p'(z_\bullet)|$ and $|p'(z_*)|$, where $p'(z) = 2z$ is the derivative of \hbox{$p(z) = z^2+\frac{i}{4}$. } 
Use the behavior of linear maps, as shown in Figure \ref{FIG_ITERATING_LINEAR_MAP}, to make a prediction about the behavior of orbits\index{orbit} for $p(z)$ near each
of these fixed points\index{fixed point}.
\end{exercise}

\begin{exercise} Let $z_\bullet$ be the fixed point\index{fixed point} for $p(z)$ discovered in Exercise
\ref{EX_FIXED_POINTS}.  Prove that for any point $z_0$ sufficiently close to
$z_\bullet$ the orbit\index{orbit} $\{z_n\}$ under $p(z) = z^2 +\frac{i}{4}$ converges to $z_\bullet$.
(I.e., prove that there exists $\delta > 0$ such that for any $z_0$ satisfying $|z_0 -z_\bullet| < \delta$ and any $\epsilon > 0$
there exists $N \in \mathbb{N}$ such that for all $n \geq N$ we have $|z_n - z_\bullet| < \epsilon$.)

Why does your proof fail if you replace the fixed point\index{fixed point} $z_\bullet$ with $z_*$?

Now, prove that the orbit\index{orbit} of $z_0 = i$ converges to $z_\bullet$.
\end{exercise}

\begin{exercise}  
Prove that there exists $r > 0$ such that for any initial condition\index{initial condition} $z_0$ with $|z_0| > r$ the orbit\index{orbit} $\{z_n\}$ of $z_0$ under $p(z) = z^2 +\frac{i}{4}$
escapes to infinity.  (I.e., prove that there exists $r > 0$ such that for any $z_0$ satisfying $|z_0| >  r$ and any $R > 0$
there exists $N \in \mathbb{N}$ such that for all $n \geq N$ we have $|z_n| > R$.)

Now prove that the orbit\index{orbit} of $w_0 = 1.1i$ escapes to infinity.
\end{exercise}

\end{example}

\begin{example}{\bf Exploring the dynamics of $\bf p_{c}: \C \rightarrow \C$ for $\bf c=-1$.}\label{EXAMPLE_C_MINUS_ONE}
In Figure \ref{FIG_ZSQUARED_MINUS1} we show the first few iterates under $p(z) = z^2-1$ of
two different orbits\index{orbit}: $\{z_n\}$ of initial condition\index{initial condition} $z_0 \approx 0.08+0.66i$ and $\{w_n\}$ of initial condition\index{initial condition} $w_0 = \frac{\sqrt{2}}{2}(1+i)$.
Orbit $\{z_n\}$
seems to converge to a periodic behavior (`periodic orbit\index{orbit}') while $\{w_n\}$ seems to escape to
$\infty$.

\begin{figure}[h!]
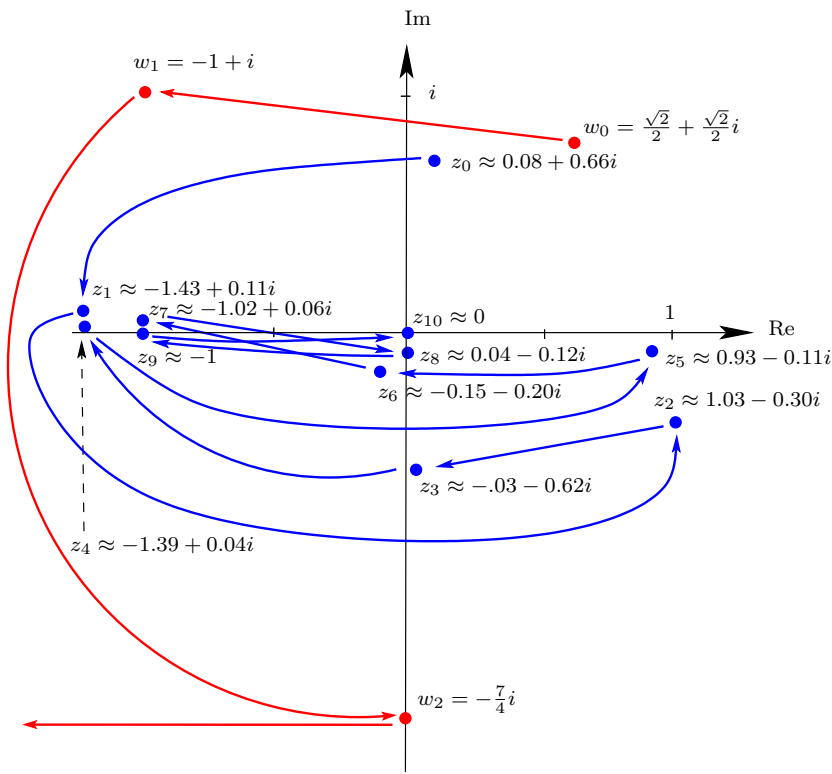
\caption{Orbits $\{z_n\}$ of $z_0 \approx 0.08+0.66i$ and $\{w_n\}$ of $w_0 = \frac{\sqrt{2}}{2}(1+i)$ under the quadratic
polynomial $p_{-1}(z) = z^2-1$.
\label{FIG_ZSQUARED_MINUS1}}
\end{figure}

In fact, the periodic orbit\index{orbit} that $\{z_n\}$ seems to converge to is easy to find for this mapping.
If we use initial condition\index{initial condition} $u_0 = 0$ we 
have $$u_1 = p_{-1}(u_0) = 0^2-1=-1.$$  Then,
\begin{align*}
u_2 = p(u_1) = p(-1) = (-1)^2-1 = 0 = u_0.
\end{align*}
  We conclude that the orbit\index{orbit} of $u_0 = 0$
is periodic with period two:
\begin{align*}
\xymatrix{
0  \ar[r]^{p_{-1}} & -1 \ar@(dl,dr)[l]^{p_{-1}}
}
\end{align*}
(Subsequently, this periodic orbit\index{orbit} will be denoted $0 \leftrightarrow 1$.)

The following two exercises are in the context of Example \ref{EXAMPLE_C_MINUS_ONE}.

\begin{exercise}  Make precise the statement that if $z_0$ is an initial condition\index{initial condition} sufficiently close to $0$, then its orbit\index{orbit} ``converges to
the periodic orbit\index{orbit} $0 \leftrightarrow 1$''.  Prove the statement.

Now, suppose $z_0 \approx 0.08+0.66i$ and prove that its orbit\index{orbit} converges to the periodic orbit\index{orbit} $0 \leftrightarrow 1$.
\end{exercise}

\begin{exercise} Find an initial condition\index{initial condition} $z_0 \in \C$ such that for any $\epsilon > 0$
there are 
\begin{enumerate}
\item infinitely many initial conditions\index{initial condition} $w_0$ with $|w_0 - z_0| < \epsilon$ having orbit\index{orbit} $\{w_n\}$ under $p_{-1}$ that remains bounded, and 
\item infinitely many initial conditions\index{initial condition} $u_0$ with $|u_0 -z_0| < \epsilon$ having orbit\index{orbit} $\{u_n\}$ under $p_{-1}$ that escapes to infinity.
\end{enumerate}
Hint: work within $\mathbb{R}$ and consider the graph of $p(x) = x^2-1$.
\end{exercise}
\end{example}

\begin{example}{\bf Exploring the dynamics of $\bf p_{c}: \C \rightarrow \C$ for $\bf c = \frac{1}{2}$.}\label{EXAMPLE_C_ONE_HALF}
As in the previous two examples, we will try a couple of arbitrary initial conditions\index{initial condition}.
Figure \ref{FIG_ZSQUARED_PLUSONEHALF} shows the orbits\index{orbit} of initial conditions\index{initial condition} $z_0 = 0$
and $w_0 \approx 0.4+0.6i$ under $p(z) = z^2+\frac{1}{2}$.  Both orbits\index{orbit} seem to escape to $\infty$.

\begin{figure}[h!]
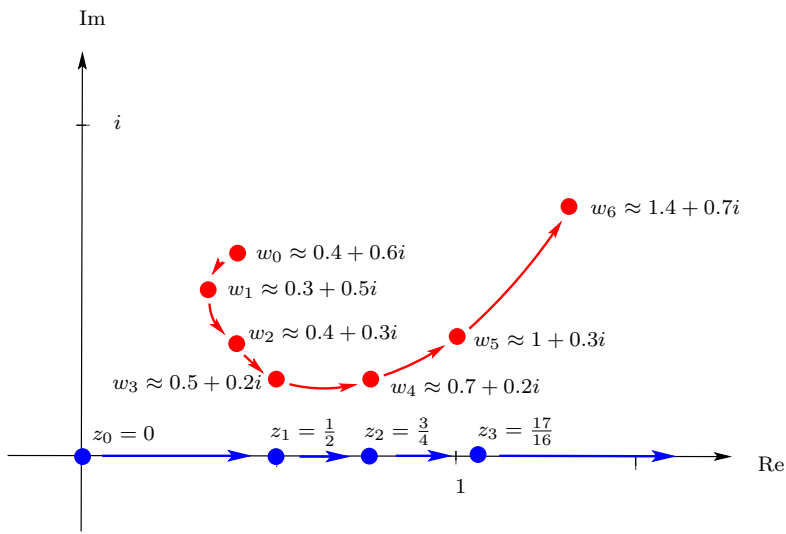
\caption{Orbits $\{z_n\}$ of $z_0 = 0$ and $\{w_n\}$ of $w_0 \approx 0.4+0.6i$ under the quadratic
polynomial $p_{\frac{1}{2}}(z) = z^2+\frac{1}{2}$.
\label{FIG_ZSQUARED_PLUSONEHALF}}
\end{figure}

\begin{exercise}\label{EX_REAL_ESCAPE} Prove that for any real initial condition\index{initial condition} $z_0 \in \mathbb{R}$ the orbit\index{orbit} $\{z_n\}$ under $p(z) = z^2+\frac{1}{2}$ escapes to infinity.
\end{exercise}

\begin{exercise}
Determine whether there is any initial condition\index{initial condition} $z_0$ for which the orbit\index{orbit} under $p_{1/2}$ remains bounded.
\end{exercise}

\end{example}

\begin{exercise}\label{EX_EXPLORE} Repeat the type of exploration done in Examples \ref{EXAMPLE_C_I_OVER_4} - \ref{EXAMPLE_C_ONE_HALF} for 
\begin{align*}
c=0, \quad c= -2, \quad c=i, \quad \mbox{and} \quad c =  -0.1 + 0.75i.
\end{align*}
Try other values of $c$.
\end{exercise}

\subsection{Questions}\label{SEC_QUESTIONS}

During our explorations we've discovered several questions.  Some of them were answered in the exercises, but several of them are still open, including:
\begin{enumerate}
\item Does every quadratic map have some initial condition\index{initial condition} $z_0$ whose orbit\index{orbit} escapes to $\infty$?

\item Does every quadratic map have some periodic orbit\index{orbit} $$z_0 \mapsto z_1
\mapsto z_2 \mapsto \cdots \mapsto z_n \mapsto z_0$$ which attracts the orbits\index{orbit}
of nearby initial conditions\index{initial condition}?  \Big(Perhaps we didn't look hard enough for one when 
$c=\frac{1}{2}$?\Big) 

\item Can a map $p_c(z)$ have more than one such attracting\index{periodic orbit!attracting} periodic orbit?\index{orbit}

\item For any $m \geq 1$ does there exist a parameter\index{parameter} $c$ such that $p_c(z)$ has an attracting\index{periodic orbit!attracting} periodic orbit of
period $m$?

\end{enumerate}

\begin{exercise}\label{EX_ESCAPE} Answer Question 1 by showing that for any $c$ there is a radius $R(c)$
such that for any initial condition\index{initial condition} $z_0$ with $|z_0| > R(c)$ the orbit\index{orbit} $\{z_n\}$ escapes to $\infty$.

Generalize your result to prove that for any polynomial $q(z)$  of degree at least $2$ there is
some $R > 0$ so that any initial condition\index{initial condition} $z_0$ with $|z_0| > R$ has orbit\index{orbit} $\{z_n\}$ that escapes to $\infty$.
\end{exercise}

\subsection{Crash course in complex analysis}

In order to answer the questions posed in the previous subsection and explore
the material more deeply, we will need some basic tools from complex analysis.
We have slightly adapted the follows results from the textbook by Saff and
Snider \cite{SS}.  We present at most sketches of the proofs and leave many of
the details to the reader.

This subsection is rather terse.  The reader may want to
initially skim over it and then move forward to see how the material is used in the later
lectures.  

\vspace{0.1in}

We begin with some topological properties of $\C$.  
The {\em open disc} of radius $r > 0$ centered at $z_0$ is 
$D(z_0,r) := \{z \in \C \,:\, |z-z_0| < r\}$.
\begin{definition}
A set $S \subset \C$ is {\em open} if for every $z \in S$ there exists $r > 0$ such that $D(z,r) \subset S$.
A set $S \subset \C$ is {\em closed} if its complement $\C \setminus S$ is open.
\end{definition}

\begin{exercise}
Prove that for any $z_0 \in \C$ and any $r > 0$ the ``open disc'' $D(z_0,r)$ is {\em actually} open.  Then prove
that the set
$$\overline{D(z_0,r)} := \{z \in \C \,:\, |z-z_0| \leq r\}$$
is closed.  It is called the {\em closed disc\index{closed disc}} of radius $r$ centered at~$z_0$.
\end{exercise}

\begin{definition}\index{boundary}
The {\em boundary} of $S \subset \C$ is 
\begin{align*}
\partial S :=  \{z \in \C \, : \, \mbox{$D(z,r)$ contains points in $S$ and in $\C \setminus S$ for every $r >0$}\}.
\end{align*}
\end{definition}

\begin{definition}
A set $S \subset \C$ is {\em disconnected\index{disconnected}} if there exist open sets $U$ and $V$ with 
\begin{itemize}
\item[(i)] $S \subset U \cup V$,
\item[(ii)] $S \cap U \neq \emptyset$ and $S \cap V \neq \emptyset$, and
\item[(iii)] $U \cap V = \emptyset$.
\end{itemize}
A set $S \subset \C$ is {\em connected\index{connected}} if it is not disconnected\index{disconnected}.
\end{definition}

An open connected\index{connected} $U \subset \mathbb{C}$ is called a {\em domain\index{domain}}.  Any set denoted $U$ in this subsection will be assumed to be a domain\index{domain}.
If $z_0 \in U$, a {\em neighborhood\index{neighborhood}} of $z_0$ will be another domain\index{domain} $V \subset
U$ with $z_0 \in V$.  (A round disc $D(z_0,r)$ for some $r > 0$ sufficiently
small will always suffice.)

\begin{definition}  A {\em contour\index{contour}} $\gamma \subset U$ is a piecewise smooth function \hbox{$\gamma: [0,1] \rightarrow U$}.  (Here, the notation implicitly
identifies the function \hbox{$\gamma: [0,1] \rightarrow U$} with its image $\gamma[0,1] \equiv \gamma \subset U$.)

A contour\index{contour} $\gamma$ is {\em closed\index{contour!closed}} if
$\gamma(0) = \gamma(1)$.  A closed contour $\gamma$ is {\em
simple\index{contour!simple}} if $\gamma(s) \neq \gamma(t)$ for $t \neq s$
unless $t=0$ and $s=1$ or vice-versa.  (Informally, a simple closed
contour as a loop that does not cross itself.) 

 A simple closed contour is {\em
positively oriented\index{contour!positively oriented}} if as you follow the
contour, the region it encloses is on your left.  (Informally, this means that it goes counterclockwise.)
\end{definition}

\begin{remark} An open set $S$ is connected\index{connected} if and only if for every two points $z,w \in S$ there is a contour\index{contour} $\gamma \subset S$ with $\gamma(0) = z$ and $\gamma(1) = w$.
\end{remark}

\begin{definition}\label{DEF_SC} A domain\index{domain} $U \subset \C$ is {\em simply connected\index{domain!simply connected}} if any closed contour\index{contour} $\gamma \subset U$ can be continuously deformed within $U$ to some point $z_0 \in U$.
\end{definition}

\noindent
We refer the reader to \cite[Section 4.4, Definition 5]{SS} for the formal definition of {\em continuously deformed\index{contour!continuously deformed}}.  In these notes, we will only need that the disc $D(z_0,r)$ is simply connected.  It follows from
the fact that any closed contour $\gamma \subset D(z_0,r)$ can be affinely scaled within $D(z_0,r)$ down to the center $z_0$.

\begin{remark} You have seen Definition \ref{DEF_SC} in your multivariable calculus class, where it was used in the statement of Green's Theorem.
\end{remark}

\begin{definition}
A set $K \subset \C$ is {\em compact\index{compact}} if for any collection
$\{W_\lambda\}_{\lambda \in \Lambda}$ of open sets with $$K \subset
\bigcup_{\lambda \in \Lambda} W_\lambda$$ there are a finite number of sets
$W_{\lambda_1},\ldots,W_{\lambda_n}$ so that $$K \subset W_{\lambda_1} \cup
\cdots \cup W_{\lambda_n}.$$  
\end{definition}

\begin{HBtheorem}
A set $S \subset \mathbb{C}$ is compact\index{compact} if and only if it is closed and bounded.
\end{HBtheorem}

\begin{exercise}\label{EX_NESTED} Suppose $K_1 \supset K_2 \supset K_3 \supset \cdots$ is a nested sequence of non-empty connected\index{connected} compact\index{compact} sets in $\C$.  Prove that $\bigcap_{n \geq 1} K_n$ is non-empty and connected\index{connected}.
\end{exercise}

\vspace{0.1in}

We are now ready to start doing complex calculus.
The notion of limit is defined in exactly the same as in calculus, except that modulus $| \cdot |$ takes the place of absolute value.

\begin{definition}
Let $z_0 \in U$ and let $f: U \setminus \{z_0\} \rightarrow \mathbb{C}$ be a function.  We say that $\lim_{z \rightarrow z_0} f(z) = L$ for some $L \in \mathbb{C}$ if for every $\epsilon > 0$ there is a $\delta > 0$ such that $0 < |z-z_0| < \delta$ implies $|f(z) -L| < \epsilon$.
\end{definition}

If we write
\begin{align*}
f(z) = u(x,y) + i v(x,y)
\end{align*}
with $u: \mathbb{R}^2 \rightarrow \mathbb{R}$ and $v: \mathbb{R}^2 \rightarrow \mathbb{R}$, then $\lim_{z \rightarrow z_0}f(z) = L$ if and only if
\begin{align*}
\lim_{(x,y) \rightarrow (x_0,y_0)} u(x,y) = \Re(L) \quad \mbox{and} \quad \lim_{(x,y) \rightarrow (x_0,y_0)} v(x,y) = \Im(L).
\end{align*}
(The limits on the right hand side are taken as in the sense of your multivariable calculus class.) 

\begin{definition}
$f: U \rightarrow \mathbb{C}$ is {\em continuous\index{continuous}} if for every $z_0 \in U$ we have $\lim_{z \rightarrow z_0} f(z) = f(z_0)$.
\end{definition}

\begin{definition}\label{DEF_DIFF}
$f: U \rightarrow \mathbb{C}$ is {\em differentiable\index{differentiable}} at $z_0 \in U$ if 
\begin{align*}
f'(z_0) := \lim_{h \rightarrow 0} \frac{f(z_0+h) - f(z_0)}{h}
\end{align*}
exists.  
\end{definition}

\begin{remark} The usual rules for differentiating sums, products, and quotients, as well
as the chain rule hold for complex derivatives.  They are proved in the same way as in your calculus class.
\end{remark}

\begin{remark} It is crucial in Definition \ref{DEF_DIFF} that one allows $h$ to approach $0$ from any direction and that the resulting limit is independent of that direction.
\end{remark}

Now for the most important definition in this whole set of notes:

\begin{definition}\label{DEF_AN}
$f: U \rightarrow \mathbb{C}$ is {\em analytic\index{analytic}} (or {\em holomorphic\index{holomorphic}}) if it is differentiable\index{differentiable} at every $z_0 \in U$.
\end{definition}

We will see that analytic\index{analytic} functions have marvelous properties!  It will be the reason why studying the iteration of analytic\index{analytic} functions is so fruitful.

\begin{exercise}
Show that $f(z) = z$ is analytic\index{analytic} on all of $\C$ and that $g(z) = \bar{z}$ is not analytic\index{analytic} in a neighborhood\index{neighborhood} of any point of $\C$.  (In fact, it is ``anti-analytic\index{anti-analytic}''.)
\end{exercise}

\begin{exercise}\label{EX_POLYS_ANALYTIC}
Show that any complex polynomial\index{complex polynomial} $$p(z) = a_d z^d + a_{d-1} z^{d-1} + \cdots + a_1 z + a_0$$ gives an analytic\index{analytic} function $p: \mathbb{C} \rightarrow \mathbb{C}$.
\end{exercise}

\begin{definition} Suppose $U$ and $V$ are domains\index{domain}.  A mapping $f: U \rightarrow V$ is called {\em conformal\index{conformal}} if it is analytic\index{analytic} and has an analytic\index{analytic} inverse $f^{-1}: V \rightarrow U$.
\end{definition}

\begin{CRtheorem}
Let $f: U \rightarrow \mathbb{C}$ be given by $$f(z) = u(x,y)+iv(x,y)$$ with $\frac{\partial u}{\partial x}, \frac{\partial u}{\partial y}, \frac{\partial v}{\partial x},$ and $\frac{\partial v}{\partial y}$ continuous\index{continuous} on $U$.  Then
\begin{align*}
f \quad \mbox{is analytic on $U$} \quad \Leftrightarrow \quad \frac{\partial u}{\partial x} = \frac{\partial v}{\partial y} \,\, \mbox{and} \,\, \frac{\partial u}{\partial y} = -\frac{\partial v}{\partial x} \,\, \mbox{for all} \quad (x,y) \in U.
\end{align*}
\end{CRtheorem}

\begin{IFtheorem}
Suppose $f: U \rightarrow \mathbb{C}$ is analytic\index{analytic} and $f'(z_0)~\neq~0$.  Then,
there is an open neighborhood\index{neighborhood} $V$ of $f(z_0)$ in $\mathbb{C}$ and an analytic\index{analytic}
function $g: V \rightarrow U$ such that $g(f(z_0)) = z_0$ and for all $w \in V$ we have $f(g(w)) = w$
and all $z \in g(V)$ we have $g(f(z)) = z$.  Moreover,
\begin{align*}
g'(f(z_0)) = \frac{1}{f'(z_0)}.
\end{align*}
\end{IFtheorem}

\begin{exercise} Show that $f(z) = z^2-1$ satisfies the hypotheses of the inverse function theorem for any $z \neq 0$.  Use the quadratic equation to
explicitly find the function $g(z)$ whose existence is asserted by the Inverse Function Theorem\index{Inverse Function Theorem}.  What goes wrong with $g$ at $-1 = f(0)$? 
\end{exercise}

\begin{EandL} According to Euler's Formula, if $z = x+iy$ with $x,y \in \mathbb{R}$ then 
\begin{align*}
\e^z = \e^{x+iy} = e^x \left(\cos y + i \sin y\right),
\end{align*}
which can be verified to be analytic\index{analytic} on all of $\C$ by using the Cauchy-Riemann Equations\index{Cauchy-Riemann Equations}.
It satisfies $(e^z)' = e^z$, which is never $0$.

Let
$S:=\{z \in \C \, : \, -\pi < \Im(z) < \pi\}$ and $\C^\dagger := \C \setminus (-\infty,0]$.
Then, the exponential function maps the strip $S$ bijectively onto $\C^\dagger$.
Therefore, it has an inverse function
\begin{align*}
{\rm Log}(z) : \C^\dagger \rightarrow S,
\end{align*}
which is analytic\index{analytic} by the inverse function theorem.  (This function is called the {\em Principal Branch of the Logarithm}.
One can define other branches that are analytic\index{analytic} on domains\index{domain} other than $\C^\dagger$; see \cite[Section 3.3]{SS}.)
\end{EandL}

\begin{definition} Suppose $f: U \rightarrow \C$ is an analytic\index{analytic} function.  A point $z \in U$
with $f'(z) = 0$ is called a {\em critical point\index{critical point}} of $f$.  A point $w \in \C$ with $w = f(z)$
for some critical point\index{critical point} $z$ is called a {\em critical value\index{critical value}}.
\end{definition}

The neighborhood\index{neighborhood} $V$ provided by the inverse function theorem could be very small.  When
combined with the Monodromy Theorem \cite[p. 295-297]{AHLFORS}, one can control the size of the domain\index{domain},
so long as it is simply connected\index{domain!simply connected}:

\begin{SCIFtheorem}
Suppose $f: U \rightarrow \C$ is an analytic\index{analytic} function  and $V \subset f(U)$ is
a simply connected\index{domain!simply connected} domain\index{domain} that doesn't contain any of the critical values\index{critical value} of
$f$. 

 Given any $w_\bullet \in V$ and any $z_\bullet \in f^{-1}(w_\bullet)$ there is a unique analytic\index{analytic}
function $g: V \rightarrow \C$  with $g(w_\bullet) = z_\bullet$, $f(g(w)) = w$ for all $w \in V$, and $g(f(z)) = z$ for all
$z \in g(V)$.
\end{SCIFtheorem}

\begin{remark} Our name for the previous result is not
standard.  Use it with caution!
\end{remark}

If $f: U \rightarrow \C$ is continuous\index{continuous} and $\gamma \subset U$ is a contour\index{contour}, then the integral
\begin{align*}
\int_\gamma f(z) dz
\end{align*}
is defined in terms of a suitable complex version of Riemann sums; see \cite[Section 4.2]{SS}.  For our purposes, we can take as definition
\begin{align*}
\int_\gamma f(z) dz := \int_{0}^1 f(\gamma(t)) \gamma'(t) dt,
\end{align*}
which is stated as Theorem 4 from \cite[Section 4.2]{SS}.

\begin{exercise} Let $\gamma$ be the positively oriented unit circle in $\mathbb{C}$.  Show that
\begin{align}\label{EQN_GOOD_CONTOUR_INTEGRAL}
\int_\gamma \frac{dz}{z} = 2\pi i,
\end{align}
which is perhaps ``the most important contour\index{contour} integral''.
\end{exercise}

\begin{Ctheorem}
If $f: U \rightarrow \C$ is analytic\index{analytic} and $U$ is simply connected\index{domain!simply connected}, then for any closed contour\index{contour} $\gamma \subset D$ we have
\begin{align*}
\int_\gamma f(z)dz = 0.
\end{align*}
\end{Ctheorem}

\begin{proof}[Sketch of proof:] The following is ``cribbed'' directly from \cite[p. 192-193]{SS}.
Write 
\begin{align*}
f(z) = u(x,y) + i v(x,y)  \qquad \mbox{and} \qquad \gamma(t) = (x(t),y(t)). 
\end{align*}
Then,
\begin{align*}
\int_\gamma f(z) dz &= \int_0^1 f(\gamma(t)) \gamma'(t) dt \\ &= \int_0^1 \Big(u(x(t),y(t))+i v(x(t),y(t))\Big)\left(\frac{dx}{dt}+ i \frac{dy}{dt}\right)dt \\
&= \int_0^1 \Big(u(x(t),y(t)) \frac{dx}{dt} - v(x(t),y(t)) \frac{dy}{dt}\Big)dt \\ & \,\,\, + i \int_0^1 \Big(v(x(t),y(t)) \frac{dx}{dt} + u(x(t),y(t)) \frac{dy}{dt}\Big)dt. 
\end{align*}
The real and imaginary parts above are just the parameterized versions of the real contour\index{contour} integrals
\begin{align*}
\int_\gamma u(x,y) dx - v(x,y) dy \qquad \mbox{and} \qquad \int_\gamma v(x,y) dx + u(x,y) dy
\end{align*}
considered in your multivariable calculus class.  Since $U$ is simply connected\index{domain!simply connected}, Green's Theorem \cite{STEWART} gives
\begin{align*}
\int_\gamma u(x,y) dx - v(x,y) dy &= \int \int_D \left(-\frac{\partial v}{\partial x} - \frac{\partial u}{\partial y}\right) dxdy \quad \mbox{and} \\
\int_\gamma v(x,y) dx + u(x,y) dy &= \int \int_D \left(\frac{\partial u}{\partial x} - \frac{\partial v}{\partial y} \right)dxdy.
\end{align*}
Since $f$ is analytic\index{analytic}, the Cauchy-Riemann Equations\index{Cauchy-Riemann Equations} imply that both integrands are $0$.  Thus, $\int_\gamma f(z)dz = 0$.
\end{proof}

\begin{remark} In the proof we have used the additional assumption that
the partial derivatives $\frac{\partial u}{\partial x}, \frac{\partial
u}{\partial y}, \frac{\partial v}{\partial x},$ and $\frac{\partial v}{\partial
y}$ are all continuous\index{continuous} functions of $(x,y)$.  This was needed in order for us
to apply Green's Theorem.  This hypothesis is not needed, but the general proof of Cauchy's Theorem\index{Cauchy's Theorem} is more complicated; see, for example, \cite[Section 4.4]{AHLFORS}.
\end{remark}

There is also an amazing ``converse'' to Cauchy's Theorem\index{Cauchy's Theorem}

\begin{Mtheorem}
If $f: U \rightarrow \C$ is continuous\index{continuous} and if  
\begin{align*}
\int_\gamma f(z)dz = 0.
\end{align*}
for any closed contour\index{contour} $\gamma \subset U$, then $f$ is analytic\index{analytic} in $U$.
\end{Mtheorem}

\begin{CIFtheorem}
Let $\gamma$ be a simple closed positively oriented contour\index{contour}.  If $f$ is analytic\index{analytic} in some simply connected\index{domain!simply connected} domain $U$ containing $\gamma$ and $z_0$ is any
point inside of $\gamma$, then
\begin{align*}
f(z_0) = \frac{1}{2\pi i} \int_\gamma \frac{f(z)}{z-z_0}dz
\end{align*}
\end{CIFtheorem}

\begin{proof}[Sketch of proof:]
Refer to Figure \ref{FIG_CIF} throughout the proof.   For any $\epsilon > 0$ we can apply Cauchy's Theorem\index{Cauchy's Theorem} to the contour\index{contour} $\eta$ proving that
\begin{align*}
\int_\gamma \frac{f(z)}{z-z_0}dz  = \int_{\gamma'} \frac{f(z)}{z-z_0}dz,
\end{align*}
where $\gamma'$ is the positively oriented circle $|z-z_0| = \epsilon$.
Since $f(z)$ is analytic\index{analytic} it is continuous\index{continuous}, implying that if we choose $\epsilon > 0$ sufficiently small, \hbox{$f(z) \approx f(z_0)$} on $\gamma'$.
Then,
\begin{align*} 
\int_{\gamma'} \frac{f(z)}{z-z_0}dz \approx f(z_0) \int_{\gamma'} \frac{1}{z-z_0} dz = 2\pi i f(z_0),
\end{align*}
with the last equality coming from (\ref{EQN_GOOD_CONTOUR_INTEGRAL}).

\begin{figure}[h!]
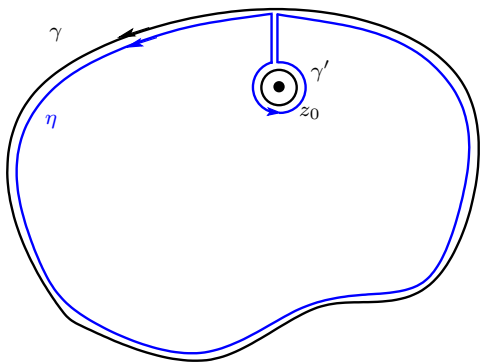
\caption{Illustration of the proof of the Cauchy Integral Formula\index{Cauchy Integral Formula}.  \label{FIG_CIF}}
\end{figure}

\end{proof}

\begin{exercise}
Use the fact that if $|f(z) - g(z)| < \epsilon$ for all $z$ on a contour\index{contour} $\gamma$ then 
\begin{eqnarray*}
\left| \int_\gamma f(z) dz - \int_\gamma g(z) dz \right| < \epsilon \, {\rm length}(\gamma)
\end{eqnarray*}
to make rigorous the estimates $\approx$ in the proof of the Cauchy Integral Formula\index{Cauchy Integral Formula}.
\end{exercise}

Let us write the Cauchy Integral Formula\index{Cauchy Integral Formula} slightly differently:
\begin{align}\label{EQN_CIF_V2}
f(z) = \frac{1}{2\pi i} \int_\gamma \frac{f(\zeta)}{\zeta-z}d\zeta,
\end{align}
where $z$ is any point inside of $\gamma$.  (This makes it more clear that we think of $z$ as an independent variable.)
By differentiating under the integral sign (after checking that it's allowed) we obtain:

\begin{CIF2theorem}
Let $\gamma$ be a simple closed positively oriented contour\index{contour}.  If $f$ is analytic\index{analytic} in some simply connected\index{domain!simply connected} domain $U$ containing $\gamma$ and $z$ is any
point inside of $\gamma$, then
\begin{align}\label{EQN_CIF_V3}
f^{(n)}(z) = \frac{n!}{2\pi i} \int_\gamma \frac{f(\zeta)}{(\zeta-z)^{n+1}}dz.
\end{align}
In particular, an analytic\index{analytic} function is infinitely differentiable\index{differentiable}!
\end{CIF2theorem}

\begin{CEtheorem}\label{COR_CAUCHY_ESTIMATES}  Suppose $f(z)$ is analytic\index{analytic} on a domain\index{domain} containing the disc $D(z_0,r)$ and suppose
$|f(z)| < M$ on the \index{boundary} boundary  $\partial D(z_0,r) = \{z \in \C\, : \, |z-z_0| = r\}$.  Then, for any $n \in \mathbb{N}$ we have
\begin{align*}
\left|f^{(n)}(z_0)\right| \leq \frac{n! M}{r^n}.
\end{align*}
\end{CEtheorem}

\begin{exercise} Prove the Cauchy Estimates\index{Cauchy Estimates}, supposing (\ref{EQN_CIF_V3}).
\end{exercise}

Suppose $\overline{D(z_0,r)} \subset U$ and $f:U \rightarrow \C$ is analytic\index{analytic}.  If we parameterize $\partial D(0,r)$
by $\gamma(t) = z_0 + r\e^{i t}$, then the Cauchy Integral Formula\index{Cauchy Integral Formula} becomes
\begin{align*}
f(z_0) = \frac{1}{2\pi} \int_0^{2\pi} f(z_0+r \e^{it}) dt.
\end{align*}
From this, one sees that it is impossible to have  $|f(z_0)| \geq |f(z_0+r \e^{it})|$ for all $t\in[0,2\pi]$ without the 
inequality actually being an equality for all $t$.  From this, it is straightforward to prove:

\begin{MMtheorem} Suppose $f(z)$ is analytic\index{analytic} in a domain\index{domain} $U$ and $|f(z)|$ achieves its maximum at a point $z_0 \in U$.
Then $f(z)$ is constant on $U$.

If, moreover, $\overline{U}$ is compact\index{compact} and $f$ extends continuously to
$\overline{U}$, then $f$ achieves its maximum modulus on the \index{boundary} boundary of $U$.
\end{MMtheorem}

Meanwhile, by using the geometric series to write
\begin{align*}
\frac{1}{\zeta-z} = \frac{1}{\zeta} \cdot \frac{1}{1-\frac{z}{\zeta}} = \frac{1}{\zeta} \, \sum_{n=0}^\infty \left(\frac{z}{\zeta}\right)^n,
\end{align*}
for any $\left|\frac{z}{\zeta}\right| < 1$,
the Cauchy Integral Formula\index{Cauchy Integral Formula} (\ref{EQN_CIF_V2}) implies:

\begin{PStheorem}Let $f$ be analytic\index{analytic} on a domain\index{domain} $U$ and suppose the disc $D(z_0,r)$ is contained in $U$.  Then, we can
write $f(z)$ as a power series
\begin{align*}
f(z) = \sum_{n=0}^\infty a_n (z-z_0)^n
\end{align*}
that converges on $D(z_0,r)$.
\end{PStheorem}

The {\em multiplicity} of a zero $z_0$ for an analytic\index{analytic} function $f(z)$ is defined as the order of the smallest non-zero
term in the power series expansion of $f(z)$ around $z_0$.

\begin{AP} Suppose $f: U \rightarrow \C$ is analytic\index{analytic} and $\gamma \subset U$ is a positively oriented simple closed contour\index{contour}
such that all points inside of $\gamma$ are in $U$.  Then, the number of zeros of $f$ (counted with multiplicities) is equal
to the change in $\arg(f(z))$ as $z$ traverses $\gamma$ once in the counter-clockwise direction.
\end{AP}

\begin{definition}\label{DEF_UC} Let $f_n: U \rightarrow \C$ be a sequence of functions and $f: U \rightarrow \C$ be another function. 
 Let $K \subset U$ be a compact\index{compact} set.
The sequence $\{f_n\}$ converges to $f$ {\em uniformly on $K$} if for every $\epsilon > 0$ there is a $\delta > 0$ such that for every $z \in K$
$|f_n(z) - f(z)| < \epsilon$.
\end{definition}

Note that the order of quantifiers in Definition \ref{DEF_UC} is crucial.  If $\delta$ was allowed to depend on $z$, we would have the weaker notion of {\em pointwise convergence}.

\begin{ULAtheorem}
Suppose $f_n: U \rightarrow \C$ is a sequence of analytic\index{analytic} functions and $f: U
\rightarrow \C$ is another (potentially non-analytic\index{analytic}) function. If for any compact\index{compact}
$K \subset U$ we have that $\{f_n\}$ converges uniformly to $f$ on $K$, then $f:
U \rightarrow \C$ is also analytic\index{analytic}. 

Moreover, for any $k \geq 1$, the $k$-th derivatives $f_n^{(k)}(z)$ converge uniformly 
to $f^{(k)}(z)$ on any compact\index{compact} $K \subset U$.
 \end{ULAtheorem}

\begin{proof}[Sketch of the proof:]
By restricting to a smaller domain\index{domain}, we can suppose $U$ is simply connected\index{domain!simply connected}.  For any contour\index{contour} $\gamma \subset U$, Cauchy's Theorem\index{Cauchy's Theorem} gives
$\int_\gamma f_n(z) dz = 0$.  Since the convergence is uniform on the compact\index{compact} set $\gamma \subset U$, we have
\begin{align*}
\int_\gamma f(z) dz = \int_\gamma \lim_{n \rightarrow \infty} f_n(z) dz =  \lim_{n \rightarrow \infty}  \int_\gamma f_n(z) dz = 0.
\end{align*}
Thus, Morera's Theorem gives that $f(z)$ is analytic\index{analytic}.

Convergence of the derivatives follows from the Cauchy Integral Formula For Higher Derivatives\index{Cauchy Integral Formula For Higher Derivatives}.
\end{proof}

The following exercises illustrate the power of the Uniform Limits Theorem\index{Uniform Limits Theorem}.

\begin{exercise}\label{EX_PS_AN}Suppose that for some $R > 0$ the power series 
\begin{align}\label{EQN_POWER_SERIES}
\sum_{n=0}^\infty a_n (z-z_0)^n
\end{align}
 converges for each $z \in D(z_0,R)$.  Prove that for
any $0 < r <R$ the power series converges uniformly on the closed disc\index{closed disc} $\overline{D(z_0,r)}$.
Use Exercise \ref{EX_POLYS_ANALYTIC} and the Uniform Limits Theorem\index{Uniform Limits Theorem} to conclude that power series (\ref{EQN_POWER_SERIES}) defines an analytic\index{analytic} function $f:D(z_0,R) \rightarrow \C$.
\end{exercise}

\begin{exercise}  Suppose we have a sequence of polynomials $p_n: [0,1] \rightarrow \mathbb{R}$ and that $p_n(x)$ converges uniformly on $[0,1]$ to some function $f: [0,1] \rightarrow \mathbb{R}$.  Does $f$ even have to be differentiable\index{differentiable}?
\end{exercise}

We close this section with the following famous result:

\begin{Slemma}
Let $\mathbb{D}:=D(0,1)$ be the unit disc and suppose $f: \mathbb{D} \rightarrow \mathbb{D}$ is analytic\index{analytic} with $f(0) = 0$.  Then
\begin{itemize}
\item[(a)] $|f'(0)| \leq 1$, and 
\item[(b)] $|f'(0)| = 1$ if and only if $f(z) = \e^{i\theta}z$ for some $\theta \in \mathbb{R}$.
\end{itemize}
\end{Slemma}

\begin{proof}[Sketch of the proof]
By the Existence of Power Series Theorem
we can write $f$ as a power series converging on $\mathbb{D}$:
\begin{align*}
f(z) = a_1 z + a_2 z^2 + a_3 z^3 \cdots,
\end{align*}
where the constant term is $0$ because $f(0) = 0$.  Therefore, 
\begin{align*}
F(z) := \frac{f(z)}{z} = a_1 + a_2 z + a_3 z^3 \cdots
\end{align*}
is also analytic\index{analytic} on $\mathbb{D}$, by Exercise \ref{EX_PS_AN}.  Applying the  Maximum Modulus Principle\index{Maximum Modulus Principle} to $F(z)$
we see that for any $0 < r < 1$ and any $\zeta$ satisfying $|\zeta| < r$ 
\begin{align*}
|F(\zeta)| \leq \frac{{\rm max}_{\{|z| = r\} }|f(z)|}{r}  \leq \frac{1}{r}.
\end{align*}
Since this holds for any $0 \leq r \leq 1$, we find that $|F(\zeta)| \leq 1$ for any $\zeta \in \mathbb{D}$.  Part (a) follows because $F(0) = f'(0)$.

If $|f'(0)| = 1$, then $|F(0)| = 1$, implying that $F$ attains its maximum at a
point of $\mathbb{D}$.  The Maximum Modulus Principle\index{Maximum Modulus Principle} implies that $F(z)$ is
constant, i.e. $F(z) = c$ for some $c$ with $|c|=1$.  Any such $c$ is of the form $\e^{i\theta}$ for some $\theta \in \mathbb{R}$, so by the definition of $F$,
we have $f(z) = \e^{i\theta} z$ for all $z \in \mathbb{D}$.

\end{proof}

\begin{remark}  There was nothing special about radius $1$.  If \hbox{$f: D(0,r) \rightarrow D(0,r)$} for some $r > 0$ and $f(0) = 0$, then (a) and (b) still hold.
\end{remark}

\section*{Lecture 2: ``Mandelbrot set from the inside out''}
\setcounter{section}{2}
\setcounter{subsection}{0}

We will work our way to the famous Mandelbrot set\index{Mandelbrot set} from an unusual perspective.

\subsection{Attracting periodic orbits\index{orbit}}

In Section \ref{SUBSEC_QUADRATIC_MAPS} we saw that the quadratic maps $p_c(z) = z^2 + c$ for $c = \frac{i}{4}, -1,$ and $-0.1+0.75i$ seemed to have
attracting\index{periodic orbit!attracting} periodic orbits\index{orbit} of periods $1, 2$, and $3$, respectively.  In this subsection we will make that notion precise and prove two results
about attracting\index{periodic orbit!attracting} periodic orbits\index{orbit}.  We will also see that the set of initial conditions\index{initial condition} whose orbits
converge to an attracting\index{periodic orbit!attracting} periodic orbit\index{orbit} can be phenomenally complicated.

While we are primarily interested in iterating quadratic polynomials $p_c(z) = z^2 +c$, it will also be helpful to 
consider iteration of higher degree polynomials $q(z)$.

\begin{definition}  A sequence 
\begin{align*}
z_0 \xrightarrow{q} z_1 \xrightarrow{q} z_2 \xrightarrow{q} \cdots \xrightarrow{q} z_m = z_0
\end{align*}
is called a {\em periodic orbit\index{periodic orbit} of period $m$} for $q$ if $z_n \neq z_0$ for each $1 \leq n \leq m-1$.
The members of such a periodic orbit\index{orbit} for $q$ are called {\em periodic points\index{periodic points}} of period $m$ for $q$.
A periodic point of period $1$ is called a {\em fixed point\index{fixed point}} of~$q$.
\end{definition}

If $z_0$ is a periodic point of period $m$ for $q$, then it is a fixed point
for the polynomial $s(z) = q^{\circ m}(z)$.  Meanwhile, if $z_0$ is a fixed point
for $s(z)$, then it is a periodic point of period $j$ for $q$, where $j$ divides 
$m$.  Thus, we can often reduce the study of periodic points to that of fixed points.

\begin{definition} 
A fixed point\index{fixed point} $z_*$ of $q$ is called {\em attracting\index{periodic orbit!attracting}} if
there is some $r > 0$ such that  such $q\big(D(z_*,r)\big) \subset D(z_*,r)$
and  for any initial condition\index{initial condition} $z_0 \in D(z_*,r)$ the orbit\index{orbit} $\{z_n\}$ under $q$
satisfies $\lim z_n = z_*$.

A periodic orbit\index{periodic orbit} $z_0 \rightarrow z_1 \rightarrow \cdots \rightarrow z_m = z_0$ is {\em attracting\index{periodic orbit!attracting}} if for each $n = 0,\ldots,m-1$ the point $z_n$ is an attracting fixed point\index{fixed point}
for $s(z) = q^{\circ m}(z)$.
\end{definition}

\begin{definition} The {\em multiplier\index{periodic orbit!multiplier}} of a periodic orbit $z_0 \rightarrow z_1 \rightarrow \cdots \rightarrow z_m = z_0$ is 
\begin{align*}
\lambda = q'(z_0) \cdot q'(z_1) \cdots q'(z_{m-1}).
\end{align*}
\end{definition}

Note that if $s(z) = q^{\circ m}(z)$, then the chain rule gives that 
\begin{align*}
s'(z_j) = q'(z_0) \cdot q'(z_1) \cdots q'(z_{m-1}) = \lambda \qquad \mbox{for each $0 \leq j \leq m-1$.}
\end{align*}
Thus the multiplier of the periodic orbit $z_0 \rightarrow z_1 \rightarrow \cdots \rightarrow z_m = z_0$ under $q$ is the same as the multiplier of each point $z_j$, when considered as a fixed point of $s(z)$.

The next lemma tells us that the same criterion we had in Section \ref{SUBSEC_QUADRATIC_MAPS}
 for $0$ being attracting\index{periodic orbit!attracting} under a linear map applies to fixed points\index{fixed point} of non-linear maps.

\begin{APPlemma} A periodic orbit\index{orbit} $$z_0 \rightarrow z_1 \rightarrow \cdots \rightarrow z_m = z_0$$ of $q$ is attracting\index{periodic orbit!attracting} if and only if its  multiplier\index{periodic orbit!multiplier} satisfies $|\lambda| < 1$.
\end{APPlemma}

\begin{proof}
Replacing $q$ by a suitable iterate we can suppose the periodic orbit\index{orbit} is a fixed point\index{fixed point} $z_*$ of $q$.  If $z_* \neq 0$ then we can consider the new polynomial $q(z+z_*) - z_*$ for which $0$ replaces $z_*$ as the fixed
point of interest.  (We call this a shift of coordinates \index{shift of coordinates}.) 

Suppose $0$ is an attracting\index{periodic orbit!attracting} fixed point\index{fixed point} for $q$.  Then, there exists $r > 0$
so that $q\big(D(0,r)\big) \subset D(0,r)$ and so that the orbit\index{orbit} $\{z_n\}$ of
any initial condition\index{initial condition} $z_0 \in D(0,r)$ satisfies $\lim_{n\rightarrow \infty}
z_n = 0$.  Since $q(0) = 0$, the Schwarz Lemma\index{Schwarz Lemma} implies that $|q'(0)| \leq 1$.
If $|q'(0)| = 1$, then the Schwarz Lemma\index{Schwarz Lemma} implies that $q$ is a rigid rotation
$z \mapsto \e^{i\theta} z$.  This would violate that the orbit\index{orbit} of any initial condition\index{initial condition} $z_0 \in
D(0,r)$ converges to $0$.  Therefore, $|q'(0)| < 1$.

Now, suppose $0$ is a fixed point\index{fixed point} for $q$ with  multiplier\index{periodic orbit!multiplier} $\lambda = q'(0)$ of modulus less than one. 
We will consider the case $\lambda \neq 0$, leaving the case $\lambda = 0$ as Exercise \ref{EX_SUPERATTRACTING}, below.
We have
\begin{align*}
q(z) = \lambda z + a_2 z^2 + \cdots + a_d z^d = \lambda \left(1+\frac{a_2}{\lambda} z + \cdots + \frac{a_d}{\lambda} z^{d-1}\right) z.
\end{align*}
Since $\lim_{z \rightarrow 0} 1+\frac{a_2}{\lambda} z + \cdots + \frac{a_d}{\lambda} z^{d-1} = 1$ and $|\lambda| < 1$ there exists $\epsilon > 0$ so that if $|z| < \epsilon$ then
\begin{align*}
\left|1+\frac{a_2}{\lambda} z + \cdots + \frac{a_d}{\lambda} z^{d-1}\right| <  1+\frac{1-|\lambda|}{2|\lambda|}.
\end{align*}
Thus, for any $|z| < \epsilon$ we have
\begin{align}\label{EQN_ATTRACTING}
|q(z)| = \left|\lambda \left(1+\frac{a_2}{\lambda} z + \cdots + \frac{a_d}{\lambda} z^{d-1}\right) \right| |z| \leq \frac{1+|\lambda|}{2} |z|.
\end{align}
In particular, $q\left(D(0,r)\right) \subset D(0,r)$ and (\ref{EQN_ATTRACTING})
implies that for any $z_0 \in D(0,r)$ the orbit\index{orbit} satisfies $|z_n| \leq \left(
\frac{1+|\lambda|}{2}\right)^n r \rightarrow 0$.  We conclude that $0$ is an
attracting\index{periodic orbit!attracting} fixed point\index{fixed point} for $q$.  
\end{proof}

\begin{exercise}\label{EX_PROVE_ATTRACTING}
Use the Attracting Periodic Orbit Lemma \index{Attracting Periodic Orbit Lemma} to verify that 
\begin{itemize}
\item[(a)] $z_* = \frac{1}{2}-\frac{\sqrt{1-i}}{2}$ is an attracting\index{periodic orbit!attracting} fixed point\index{fixed point} for $p(z) = z^2 +\frac{i}{4}$,
\item[(b)] $0 \leftrightarrow 1 $ is an attracting\index{periodic orbit!attracting} periodic orbit\index{orbit} of period $2$ for \hbox{$p(z) = z^2 -1$,} and
\item[(c)] If $c$ satisfies $c^3+2c^2+c+1=0$, then $0 \rightarrow c \rightarrow c^2+c \rightarrow 0$
is an attracting\index{periodic orbit!attracting} periodic orbit of period $3$.  (One of the solutions for $c$ is
the parameter\index{parameter} $c \approx -0.12+0.75i$ studied in Exercise~\ref{EX_EXPLORE}.)
\end{itemize}
\end{exercise}

\begin{exercise} Verify that there exists $r > 0$ such that for any initial condition\index{initial condition} $z_0 \in \mathbb{R}$ with $|z_0| < r$ the orbit\index{orbit} under $q(z) = z-z^3$ converges to $0$.  Why is $0$ not attracting\index{periodic orbit!attracting} as a complex fixed point\index{fixed point}?
\end{exercise}

\begin{exercise}\label{EX_SUPERATTRACTING}
Prove that if $z_*$ is a fixed point\index{fixed point} for a polynomial $q$ having  multiplier\index{periodic orbit!multiplier} $\lambda = 0$, then
$z_*$ is attracting\index{periodic orbit!attracting}.
\end{exercise}

\begin{definition} Suppose $\mathcal{O} = z_0 \rightarrow z_1 \rightarrow \cdots \rightarrow z_m = z_0$ is an attracting\index{periodic orbit!attracting} periodic orbit.  The {\em basin of attraction\index{periodic orbit!basin of attraction}} $\A(\mathcal{O})$ is 
\begin{align*}
\A(\mathcal{O}) := \{z \in \mathbb{C} \,:\, \mbox{$s^{\circ n}(z) \rightarrow z_j$ as $n \rightarrow \infty$ for some $0 \leq j \leq m-1$} \},
\end{align*}
where $s(z) = q^{\circ m}(z)$.
The {\em immediate basin\index{periodic orbit!immediate basin}} $\A_0(\mathcal{O})$ is the union of the connected\index{connected} components of $\A(\mathcal{O})$ containing the points $z_0,\ldots,z_{m-1}$.
\end{definition}

\begin{figure}[h!]
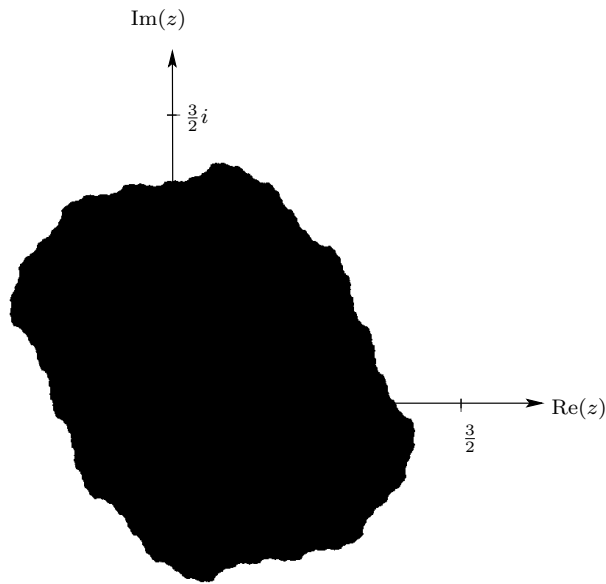
\caption{Basin of attraction of the fixed point\index{fixed point} \hbox{$z_* = \frac{1}{2}-\frac{\sqrt{1-i}}{2}$} for
$p(z) = z^2 + \frac{i}{4}$.\label{FIG_BASINS1}}
\end{figure}

\begin{figure}[h!]
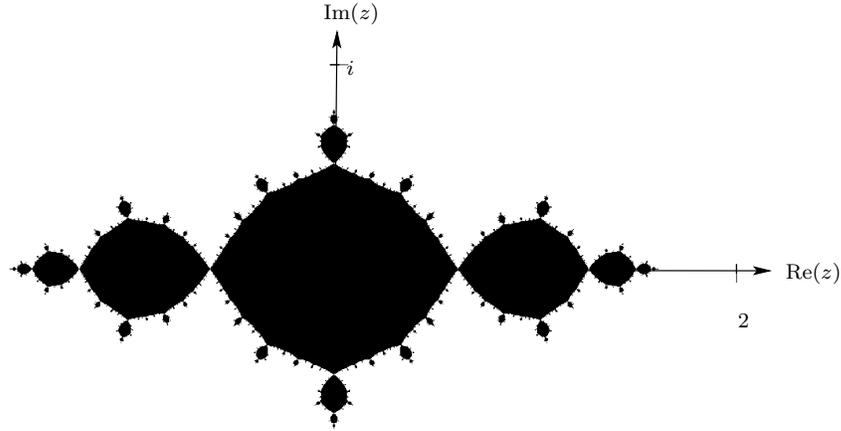
\caption{Basin of attraction of the period two cycle \hbox{$0 \leftrightarrow -1$} for
$p(z) = z^2 -1$.\label{FIG_BASINS2}}
\end{figure}

\begin{figure}[h!]
\scalebox{1.05}{
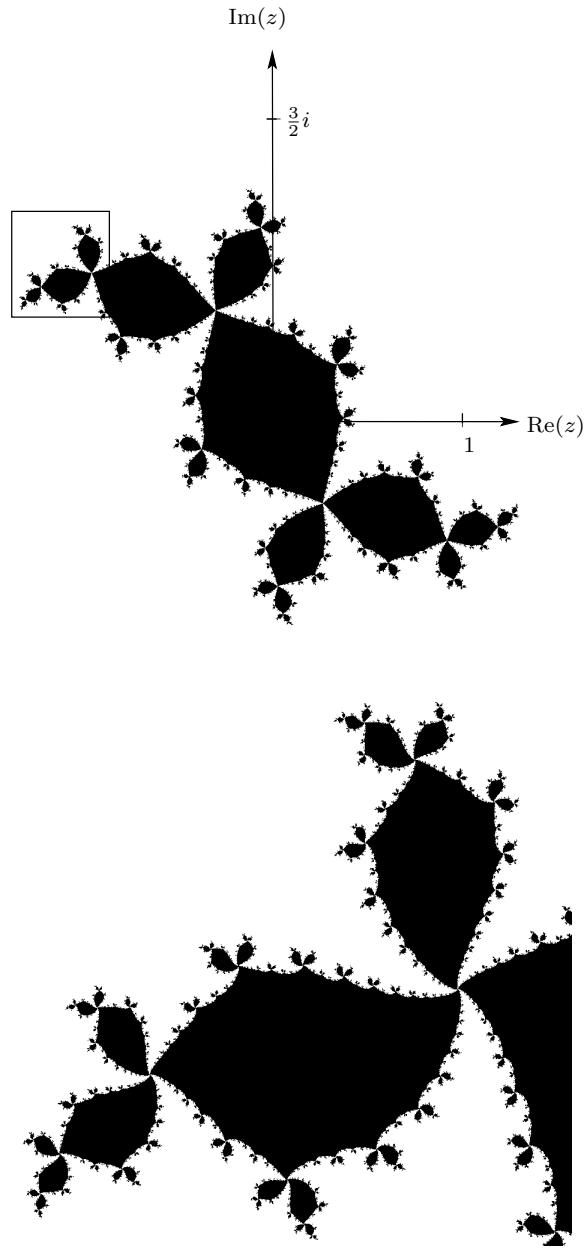
}
\caption{Top: basin of attraction of the attracting\index{periodic orbit!attracting} period $3$ cycle $0 \rightarrow c \rightarrow c^2+c$
for $c \approx -0.12+0.75i$.  Bottom: zoomed-in view of the boxed region from the left. \label{FIG_BASINS3}}
\end{figure}

Computer generated images of 
the basins of attraction for the attracting\index{periodic orbit!attracting} periodic orbits discussed in
Exercise \ref{EX_PROVE_ATTRACTING} are shown in Figures \ref{FIG_BASINS1} - \ref{FIG_BASINS3}.
Notice the remarkable complexity of the boundaries of the basins of attraction, something we would 
never have guessed during our experimentation in Section \ref{SUBSEC_QUADRATIC_MAPS}.  

\begin{figure}[h!]
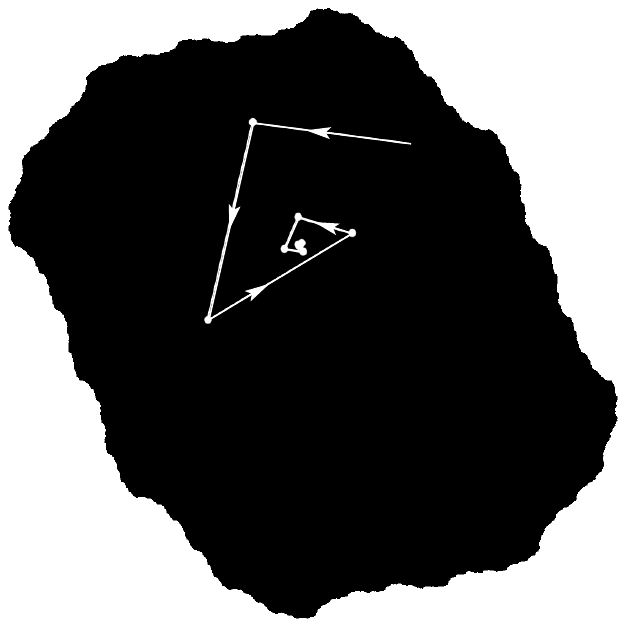
\caption{An orbit\index{orbit} converging to the attracting\index{periodic orbit!attracting} fixed point\index{fixed point} for $p(z) = z^2+\frac{i}{4}$.  Here, $\lambda = 1-\sqrt {1-2\,i} \approx 0.8 \e^{1.9i}$.  \label{FIG_NONLINEAR_SPIRAL}}
\end{figure}

\begin{figure}[h!]
\scalebox{1.4}{
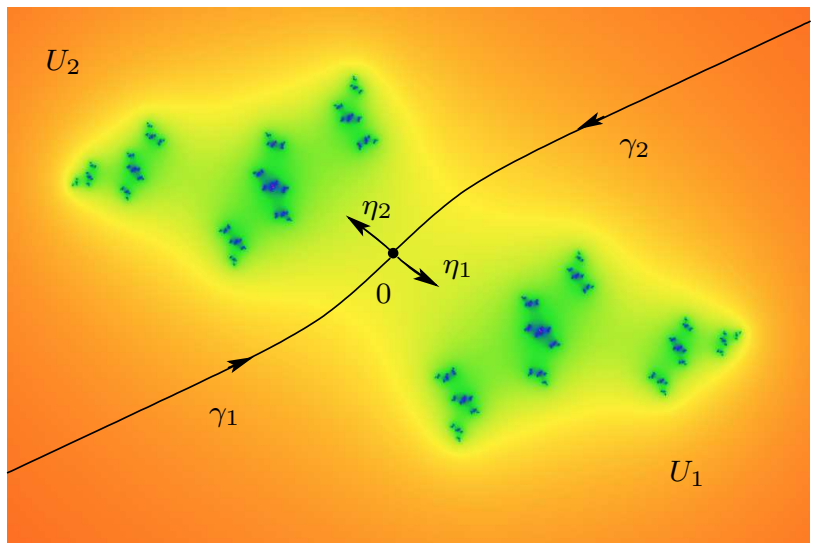
}
\caption{Stable and unstable trajectories of $-\nabla G_c$ for the critical point\index{critical point} $0$.
\label{FIG_MORSE}}
\end{figure}

\vspace{0.1in}
\noindent
{\bf Remark on computer graphics:} We used Fractalstream \cite{FS} to create Figures \ref{FIG_BASINS1}-\ref{FIG_NONLINEAR_SPIRAL}, \ref{FIG_MANDEL_BW}-\ref{FIG_MORSE}, and \ref{FIG_SMALL_MANDELBROT}-\ref{FIG_MANDEL_ZOOMED}.  Other useful programs include Dynamics Explorer \cite{DE} and the Boston University Java Applets \cite{JAVA}. \index{computer software}
\vspace{0.1in}

It is natural to ask how complicated the dynamics for iteration of $q$ can be near an attracting\index{periodic orbit!attracting}
fixed point\index{fixed point}.  The answer is provided by K\oe nig's Theorem\index{K\oe nig's Theorem} and B\"ottcher's Theorem\index{B\"ottcher's Theorem}.

\begin{koenigstheorem}\label{THM_KOENIG}
Suppose $z_\bullet$ is an attracting\index{periodic orbit!attracting} fixed point\index{fixed point} for $q$ with  multiplier\index{periodic orbit!multiplier} $\lambda
\neq 0$.  Then, there exists a neighborhood\index{neighborhood} $U$ of $z_\bullet$ and a conformal\index{conformal} map
\begin{align*}
\phi: U \rightarrow \phi(U) \subset \mathbb{C}
\end{align*}
 so that for any $w
\in \phi(U)$ we have
\begin{align}\label{EQN_PHI_DESIRED}
\phi \circ q \circ \phi^{-1} (w) = \lambda w.
\end{align}
\end{koenigstheorem}

In other words, Theorem \ref{THM_KOENIG} gives that there is a neighborhood\index{neighborhood} $U$ of
$z_\bullet$ in which there is a coordinate system $w = \phi(z)$ in which the
non-linear mapping $q$ becomes linear!  This explains why the the same
geometric spirals shown on the top of Figure
\ref{FIG_ITERATING_LINEAR_MAP} for the linear map appear sufficiently close to
an attracting\index{periodic orbit!attracting} fixed point\index{fixed point} $z_\bullet$ for a non-linear map.  This is
illustrated in Figure~\ref{FIG_NONLINEAR_SPIRAL}.

\begin{proof}
Shifting \index{shift of coordinates} the coordinates if necessary, we can suppose $z_\bullet = 0$.  
The Attracting Periodic Orbit Lemma\index{Attracting Periodic Orbit Lemma}
gives that the  multiplier\index{periodic orbit!multiplier} of $0$ satisfies $|\lambda| < 1$.  Therefore, as in
the second half of the proof of the Attracting Periodic Orbit Lemma\index{Attracting Periodic Orbit Lemma},
we can find some $r > 0$ and $|\lambda| < a < 1$ so that
\begin{align}\label{EQN_GOOD_CONVERGENCE}
\mbox{for any} \quad z \in D(0,r) \quad \mbox{we have} \quad |z_n| \leq a^n r,
\end{align}
where $z_n := q^{\circ n}(z)$.

Since $q(0) = 0$ we have
\begin{align}\label{EQN_Q_WITH_REMAINDER}
q(z) = \lambda z + s(z),
\end{align}
where $s(z) = a_2 z^2 + a_3 z^3 + \cdots a_d z^d$.  In particular, there exists $b > 0$ so that
\begin{align}\label{EQN_REMAINDER}
|s(z)| \leq b |z|^2.
\end{align}
for all $z \in D(0,r)$.

Let 
\begin{align*}
\phi_n : D(0,r) \rightarrow \C \qquad \mbox{be given by} \qquad \phi_n(z) := \frac{z_n}{\lambda^n},
\end{align*}
which satisfies $\phi_n(0) = 0$, since $0$ is a fixed point\index{fixed point}.
Notice that 
\begin{align}\label{EQN_DEF_PHIN}
\phi_n(q(z)) = \frac{z_{n+1}}{\lambda^n} =  \lambda \cdot \frac{z_{n+1}}{\lambda^{n+1}} = \lambda \phi_{n+1}(z).
\end{align}
Suppose we can prove that $\phi_n$ converges uniformly on 
$D(0,r)$ to some function $\phi: D(0,r) \rightarrow \C$.  Then,
$\phi$ will be analytic\index{analytic} by the Uniform Limits Theorem\index{Uniform Limits Theorem}.  Meanwhile, the left and right sides of~(\ref{EQN_DEF_PHIN}) converge to $\phi(q(z))$ and $\lambda \phi(z)$, respectively, implying
\begin{align}\label{EQN_PHI_DESIRED2}
\phi(q(z)) = \lambda \phi(z).
\end{align}
Since $\phi_n(0) = 0$ for each $\phi_n$ we will also have $\phi(0) = 0$.

To see that the $\phi_n$ converge uniformly on $D(0,r)$, let us rewrite it as
\begin{align*}
\phi_n(z) = \frac{z_n}{\lambda^n} = z_0 \cdot \frac{z_1}{\lambda  z_0} \cdot \frac{z_2}{\lambda  z_1} \cdot \frac{z_3}{\lambda  z_2} \cdots \frac{z_n}{\lambda  z_{n-1}}.
\end{align*}
By (\ref{EQN_Q_WITH_REMAINDER}), the general term of the product becomes 
\begin{align*}
\frac{z_k}{\lambda  z_{k-1}} = \frac{q(z_{k-1})}{\lambda  z_{k-1}} = \frac{\lambda  z_{k-1} + s(z_{k-1})}{
\lambda  z_{k-1}} = 1+ \frac{s(z_{k-1})}{\lambda z_{k-1}}.
\end{align*}
By the estimates (\ref{EQN_REMAINDER}) and (\ref{EQN_GOOD_CONVERGENCE}) on $|z_n|$  we have
\begin{align}\label{EQN_SMALL_STUFF}
\left|\frac{s(z_{k-1})}{\lambda z_{k-1}}\right| \leq b \frac{|z_k|}{\lambda} \leq b \frac{a^k r}{\lambda}.
\end{align}
We will now make $r$ smaller, if necessary, to ensure that the right hand side of (\ref{EQN_SMALL_STUFF}) is less than $\frac{1}{2}$.

To show that the $\phi_n$ converge uniformly on $D(0,r)$, it is sufficient to show that 
the infinite product
\begin{align*}
\frac{z_1}{\lambda  z_0} \cdot \frac{z_2}{\lambda  z_1} \cdot \frac{z_3}{\lambda  z_2} \cdots \frac{z_n}{\lambda  z_{n-1}} \cdots
\end{align*}
does.  Such a product converges if and only if logarithms of the finite partial products converge, i.e. if and only if the infinite sum
\begin{align}\label{EQN_LOG_SUM}
{\rm Log} \phi(z) = \sum_{k=1}^\infty {\rm Log} \left(1+\frac{s(z_{k-1})}{\lambda z_{k-1}}\right)
\end{align}
converges.  (We can take the logarithms on the right hand side of (\ref{EQN_LOG_SUM}) because
our bound of (\ref{EQN_SMALL_STUFF}) by $\frac{1}{2}$ implies that $1+\frac{s(z_{k-1})}{\lambda z_{k-1}} \in \C \setminus (-\infty,0]$.)
Using the estimate
\begin{align*}
|{\rm Log}(1+w)| \leq 2|w|  \qquad \mbox{for any} \quad |w| < \frac{1}{2}
\end{align*}
and (\ref{EQN_SMALL_STUFF})
we see that the $k$-th term is geometrically small:
\begin{align*}
\left|{\rm Log}\left(1+\frac{s(z_{k-1})}{\lambda z_{k-1}}\right)\right| \leq 2 b \frac{a^k r}{\lambda}.
\end{align*}
This proves convergence of (\ref{EQN_LOG_SUM})

It remains to show that $\phi$ is conformal\index{conformal} when restricted to a small enough neighborhood\index{neighborhood} $U \subset D(0,r)$ of $0$.
By the chain rule, each $\phi_n$ satisfies $\phi_n'(0)~=~1$.  Since the $\phi_n$ converge uniformly
to $\phi$ in a neighborhood\index{neighborhood} of $0$, the Cauchy Integral Formula\index{Cauchy Integral Formula} (\ref{EQN_CIF_V3}) implies that
$\phi_n'(0) \rightarrow \phi'(0)$.  Thus $\phi'(0) = 1$.  

By the Inverse Function Theorem\index{Inverse Function Theorem}, there is a
neighborhood\index{neighborhood} $V$ of $0 = \phi(0)$ and an
analytic\index{analytic} function $g: V \rightarrow D(0,r)$ so that $\phi(g(w)) =
w$ for every $w \in V$.  If we let $U = g(V)$, then $\phi: U \rightarrow V$ is
conformal\index{conformal}.

To obtain (\ref{EQN_PHI_DESIRED}), precompose (\ref{EQN_PHI_DESIRED2}) with $\phi^{-1} = g$ on $V$.
\end{proof}

\begin{extendedexercise}
Adapt the proof of K\oe nig's Theorem to prove:
\begin{bottcherstheorem}
Suppose $p(z)$ has a fixed point\index{fixed point} $z_\bullet$ of  multiplier\index{periodic orbit!multiplier} $\lambda = 0$ and thus is of the form
\begin{align*}
p(z) = z_\bullet + a_k(z-z_\bullet)^k + a_{k+1} (z-z_\bullet)^{k+1} + \cdots + a_d (z-z_\bullet)^d
\end{align*}
for some $2 \leq k < d$.  
Then, there exists a neighborhood\index{neighborhood} $U$ of $z_\bullet$ and a conformal\index{conformal} map
$$\phi: U \rightarrow \phi(U) \subset \mathbb{C}$$
 so that for any $w
\in \phi(U)$ we have
$\phi \circ p \circ \phi^{-1} (w) =  w^k$.
\end{bottcherstheorem}

\end{extendedexercise}

\subsection{First Exploration of the Parameter Space: The Set $M_0$}
\label{SUBSEC_FIRST_PARAM}

Let us try to understand the space of parameters\index{parameter} $c \in \C$ for the quadratic polynomial maps $p_c(z) = z^2+c$.
Consider
\begin{align*}
M_0 := \{c \in \C \, :  \, \mbox{$p_c(z)$ has an attracting\index{periodic orbit!attracting} periodic orbit} \}.
\end{align*}
We have already seen in Section \ref{SUBSEC_QUADRATIC_MAPS} that $c =
\frac{i}{4}$ and $c = -1$ are in $M_0$ and that $c = \frac{1}{2}$ is probably
not in $M_0$.
We will now use the Attracting Periodic Orbit Lemma\index{Attracting Periodic Orbit Lemma} to find some regions that are
in $M_0$. 

\begin{figure}[h!]
\scalebox{0.9}{
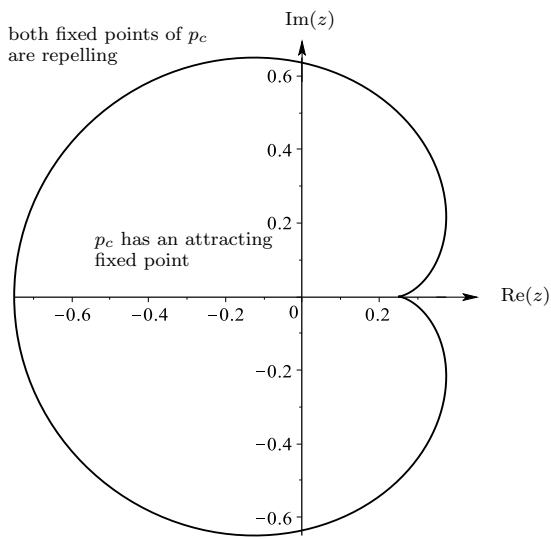
}
\caption{$p_c(z) = z^2+c$ has an attracting\index{periodic orbit!attracting} fixed point\index{fixed point} if and only if $c$ lies inside the Cardiod $c = \frac{e^{i\theta}}{2} - \frac{e^{i2\theta}}{4}$, where $0 \leq \theta \leq 2\pi$, depicted here.\label{FIG_CARDIOD}}
\end{figure}

The fixed points\index{fixed point} of $p_c(z) = z^2 + c$ are
\begin{align*}
z_* = \frac{1}{2} + \frac{\sqrt{1-4c}}{2} \qquad \mbox{and} \qquad z_\bullet = \frac{1}{2} - \frac{\sqrt{1-4c}}{2}
\end{align*}
and, since $p_c'(z) = 2z$, their  multipliers\index{periodic orbit!multiplier} are
\begin{align*}
\lambda_* = 1+\sqrt{1-4c} \qquad \mbox{and} \qquad \lambda_\bullet = 1-\sqrt{1-4c}.
\end{align*}
If $|\lambda_*| = 1$, then 
\begin{align*}
1+\sqrt{1-4c} = \e^{i\theta}
\end{align*}
for some $\theta \in \mathbb{R}$.  Solving for $c$, we find
\begin{align*}
c = \frac{e^{i\theta}}{2} - \frac{e^{i2\theta}}{4}.
\end{align*}
The resulting curve $C$ is a ``Cardiod'', shown in Figure \ref{FIG_CARDIOD}.

In each of the two regions of $\C \setminus C$ we choose the points $c=0$ and
$c=1$, which result in $\lambda_* = 2$ and $\lambda_* = 1+\sqrt{3}i$, respectively.  Thus,  neither
of the regions from $\C \setminus C$ 
corresponds to parameters\index{parameter} $c$ for which $z_*$ is an attracting\index{periodic orbit!attracting} fixed point\index{fixed point}.  Thus, we conclude that the
smallest $|\lambda_*|$ can be is $1$, occurring exactly on the Cardiod $C$. 

Doing the same computations with the  multiplier\index{periodic orbit!multiplier} $\lambda_\bullet$ of the second
fixed point\index{fixed point} $z_\bullet$,  we also find that $|\lambda_\bullet| = 1$ if and only if $c$ is on the
Cardiod $C$.  However, at $c=0$ we have $\lambda_\bullet = 0$ and at $c=1$ we
have $\lambda_\bullet = 1-\sqrt{3}i$, which is of modulus greater than $1$.
Therefore, according to the Attracting Periodic Orbit Lemma\index{Attracting Periodic Orbit Lemma}, fixed point\index{fixed point} $z_\bullet$ is
attracting\index{periodic orbit!attracting} if and only if $c$ is inside of the Cardiod $C$.  We summarize the past three paragraphs with:

\begin{lemma} $p_c(z) = z^2 + c$ has an attracting\index{periodic orbit!attracting} fixed point\index{fixed point} if and only if $c$ lies inside
the Cardiod curve $C:= \left\{c = \frac{e^{i\theta}}{2} - \frac{e^{i2\theta}}{4} \, : \, 0 \leq \theta \leq 2\pi\right\}$.
\end{lemma}

To find periodic orbits\index{periodic points} of period two, we solve 
$p^{\circ 2}_c(z) = \left( {z}^{2}+c \right) ^{2}+c = z$.   In addition to the two fixed points\index{fixed point} $z_*$ and $z_\bullet$, we find 
\begin{align*}
z_0 = -\frac{1}{2}+\frac{\sqrt {-3-4\,c}}{2} \qquad \mbox{and} \qquad z_1 =  -\frac{1}{2}-\frac{\sqrt {-3-4\,c}}{2}.
\end{align*}
One can check that $p_c(z_0) = z_1$ and $p_c(z_1) = z_0$.  These points are equal if $c = -\frac{3}{4}$, otherwise, they are indeed a periodic orbit of period $2$.

The  multiplier\index{periodic orbit!multiplier} of this periodic orbit is
\begin{align*}
\lambda = \left(-1+\sqrt {-3-4\,c}\right)\left(-1-\sqrt {-3-4\,c}\right) = 4+4c,
\end{align*}
which has modulus $1$ if and only if $|c+1| = \frac{1}{4}$.  Since $\lambda = 0$ for $c=-1$ (inside the circle) and $\lambda = 4$ for $c= 0$ (outside the circle) we find:
\begin{lemma} $p_c(z) = z^2 + c$ has an periodic orbit of period $2$ if and only if $c$ lies inside
the circle $|c+1| = \frac{1}{4}$.
\end{lemma}

\begin{figure}[h!]
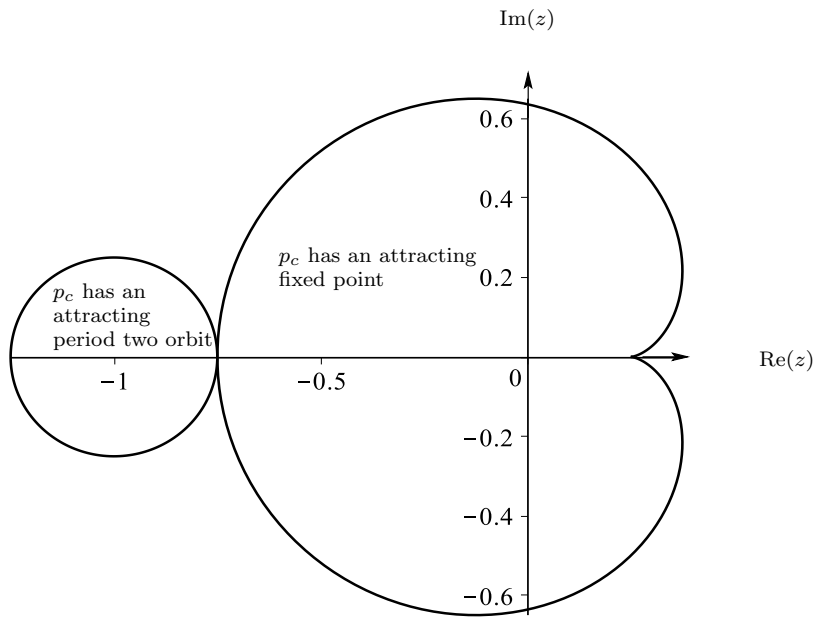
\caption{The regions in the parameter\index{parameter} plane where $p_c(z) = z^2+c$ has an attracting\index{periodic orbit!attracting} fixed point\index{fixed point} 
and where $p_c$ has an attracting\index{periodic orbit!attracting} periodic orbit of period $2$.  Combined, they form a subset of $M_0$.\label{FIG_BAD_MANDEL}}
\end{figure}

In Figure \ref{FIG_BAD_MANDEL} we show the regions of $M_0$ that we have discovered.

\begin{exercise} If possible, determine the region of parameters\index{parameter} $c$ for which $p_c(z) = z^2+c$
has an attracting\index{periodic orbit!attracting} periodic orbit of period $3$.
\end{exercise}

As $n$ increases, this approach becomes impossible.  We need a different approach, which 
requires a deeper study of attracting\index{periodic orbit!attracting} periodic orbits.

\subsection{Second Exploration of the Parameter Space: The Mandelbrot set\index{Mandelbrot set} $M$}

\begin{fatoujulialemma}
Let $q$ be a polynomial of degree $d \geq 2$.  Then, the immediate basin of
attraction for any attracting\index{periodic orbit!attracting} periodic
orbit\index{orbit} contains at least one critical point\index{critical point}
of~$q$.  In particular, since $q$ has $d-1$ critical points\index{critical
point} (counted with multiplicity), $q$ can have no more than $d-1$ distinct
attracting\index{periodic orbit!attracting} periodic orbits.
\end{fatoujulialemma}

\begin{figure}[h!]
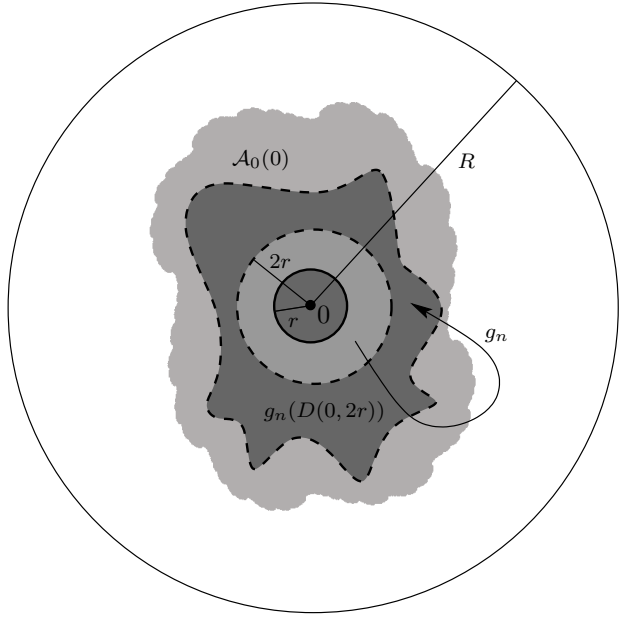
\caption{Illustration of the proof of the Fatou-Julia Lemma\index{Fatou-Julia Lemma}.
\label{FIG_FATOU_JULIA_ILLUSTRATION}}
\end{figure}

The following proof is illustrated in Figure \ref{FIG_FATOU_JULIA_ILLUSTRATION}.

\begin{proof}
Replacing $q$ by an iterate, we can suppose that the attracting\index{periodic orbit!attracting} periodic orbit is a fixed point\index{fixed point} $z_\bullet$ of $q$.  Performing a \index{shift of coordinates} shift of coordinates, we suppose $z_\bullet =0$.

If $0$ has  multiplier\index{periodic orbit!multiplier} $\lambda = 0$, then $0$ is already a critical point\index{critical point} in
the immediate basin $\A_0(0)$.  We therefore suppose $0$ has  multiplier\index{periodic orbit!multiplier}
$\lambda \neq 0$.   By the Attracting Periodic Orbit Lemma\index{Attracting Periodic Orbit Lemma}, $|\lambda| < 1$.

Suppose for contradiction that there is no critical point\index{critical point} for $q$ in $\A_0(0)$.

According to Exercise \ref{EX_ESCAPE} there is some $R > 0$ so that any initial condition\index{initial condition} $z_0$ with
$|z_0| > R$ has orbit $\{z_n\}$ that escapes to $\infty$.  In particular,
\begin{align}\label{EQN_BOUND_ON_IMMEDIATE_BASIN}
\A_0(0) \subset D(0,R).
\end{align}

We claim that $q(\A_0(0)) = \A_0(0)$. 
Since $\A(0)$ is forward invariant, $q(\A_0(0)) \subset \A(0)$.  Because $\A_0(0)$ is connected\index{connected}, so is
$q(\A_0(0))$, which is therefore contained in one of the connected\index{connected} components of $\A(0)$.
Since $0 = q(0) \in q(\A_0(0))$, we have $q(\A_0(0)) \subset \A_0(0)$.  

Conversely, suppose
$z_* \in \A_0(0)$.  Let $\gamma$ be a simple contour\index{contour} in $\A_0(0)$ connecting $z_*$ to $0$ and avoiding
any critical values\index{critical value} of $q$.  (By hypothesis, such critical values\index{critical value} would be images of critical points\index{critical point}
that are not in $\A_0(0)$.)  Then, $q^{-1}(\gamma)$ is a union of several simple contours\index{contour}.
Since $q^{-1}(0) = 0$, one of them is a simple contour\index{contour} ending at $0$.  The other end is a point
$z_\#$, which is therefore in $\A_0(0)$.  By construction $q(z_\#) = z_*$.

To simplify notation, let $f := q|_{{\A_0(0)}} : \A_0(0) \rightarrow \A_0(0)$, which satisfies
\begin{enumerate}
\item $f(\A_0(0)) = \A_0(0)$  and
\item $f$ has no critical points\index{critical point}.
\end{enumerate}

These properties persists under iteration, giving that
$f^n(\A_0(0)) = \A_0(0)$ and $f^n$ has no critical points\index{critical point} for
every $n \geq 1$.  (The latter uses the Chain Rule.)

Let $r > 0$ be chosen sufficiently small so that $D(0,2r) \subset \A_0(0)$.
Since $D(0,2r)$ is simply connected\index{domain!simply connected}, the Simply-Connected Inverse Function Theorem\index{Simply-Connected Inverse Function Theorem}
gives for each $n \geq 1$ an analytic\index{analytic} function 
\begin{align*}
g_n: D(0,2r) \rightarrow \A_0(0) \subset D(0,R)
\end{align*}
with $g_n(0) = 0$ and $f^{\circ n}(g_n(w)) = w$ for all $w \in (0,2r)$.  Its derivative satisfies
\begin{align}\label{EQN_BIG_DERIVATIVES}
g_n'(0) = \frac{1}{(f^n)'(0)} = \frac{1}{\lambda^n},
\end{align}
which can be made arbitrarily large by choosing $n$ sufficiently large, since $|\lambda| < 1$.

Meanwhile, we can apply the Cauchy Estimates\index{Cauchy Estimates} \ref{COR_CAUCHY_ESTIMATES} to the closed disc\index{closed disc} $\overline{D(0,r)} \subset D(0,2r)$.  They assert that
\begin{align*}
|g_n'(0)| \leq \frac{R}{r},
\end{align*}
where $R$ is the bound on the radius of $\A_0(0)$  given in (\ref{EQN_BOUND_ON_IMMEDIATE_BASIN}).
This is a contradiction to (\ref{EQN_BIG_DERIVATIVES}).
We conclude that the immediate basin $\A_0(0)$ contains a critical point\index{critical point} of~$q$.
\end{proof}

\begin{exercise} Use the Fatou-Julia Lemma\index{Fatou-Julia Lemma} and the result of Exercise \ref{EX_REAL_ESCAPE} to (finally) prove that $p(z) = z^2 + \frac{1}{2}$ does
not have any attracting\index{periodic orbit!attracting} periodic orbit.  This answers Question 2 from Section \ref{SEC_QUESTIONS} in the negative.
\end{exercise}

\begin{remark}
The Fatou-Julia Lemma\index{Fatou-Julia Lemma} also answers our Question 3 from Section \ref{SEC_QUESTIONS} by telling us
that a quadratic polynomial can have at most one attracting\index{periodic orbit!attracting} periodic orbit.
\end{remark}

In Section \ref{SUBSEC_FIRST_PARAM} we were interested in the set
\begin{align*}
M_0 := \{c \in \C \, : \, \mbox{$p_c(z)$ has an attracting\index{periodic orbit!attracting} periodic orbit} \}.
\end{align*}
Using the Attracting Periodic Orbit Lemma\index{Attracting Periodic Orbit Lemma} to find regions in the complex plane
for which $p_c(z) = z^2 + c$ has an attracting\index{periodic orbit!attracting} periodic point of period $n$ 
became hopeless once $n$ is large.  
The results for $n=1$ and $2$ are shown in Figure \ref{FIG_BAD_MANDEL}.

If we decide to lose control over what period the attracting\index{periodic orbit!attracting} periodic point has, the Fatou-Julia Lemma\index{Fatou-Julia Lemma} gives us some very interesting information:

\begin{corollary}{\bf (Consequence of Fatou-Julia Lemma)}\label{COR_TO_FATOU_JULIA}
If $p_c(z)$ has an attracting\index{periodic orbit!attracting} periodic orbit, then the orbit $\{p_c^{\circ n}(0)\}$ of the critical point\index{critical point} $0$ remains bounded.
\end{corollary}

\begin{figure}[h!]
\scalebox{1.3}{
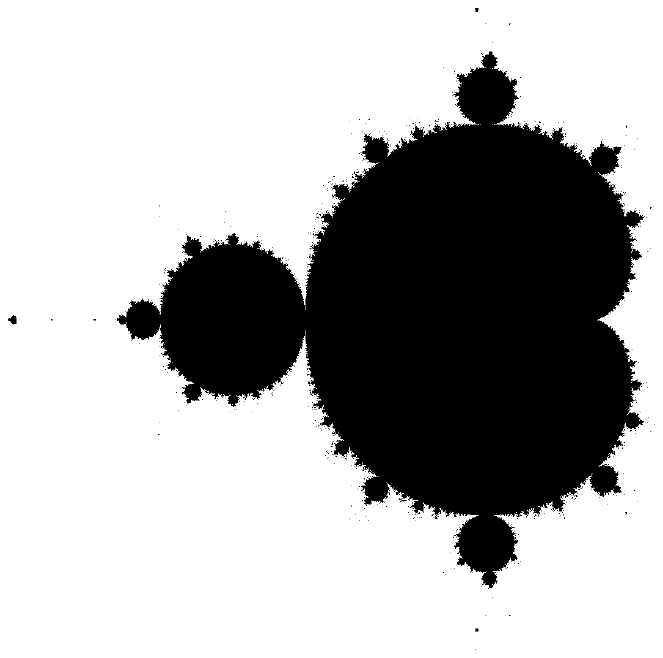
}
\caption{The Mandelbrot set\index{Mandelbrot set} $M$, shown in black.  The region shown is
approximately $-2.4 \leq \Re(z) \leq 1$ and $-1.6 \leq \Im(z) \leq 1.6$. 
\label{FIG_MANDEL_BW}}
\end{figure}

This motivates one to define another set:

\begin{definition} The {\em Mandelbrot set\index{Mandelbrot set}} is
\begin{align}\label{EQN_DEF_MANDELBROT}
M := \{c \in \C \, : \, \mbox{$p_c^{\circ n}(0)$ remains bounded for all $n \geq 0$} \}.
\end{align}
\end{definition}
\noindent
A computer image of the Mandelbrot set\index{Mandelbrot set} is depicted in Figure \ref{FIG_MANDEL_BW}.
One sees small ``dots'' at the left end and top and bottom of the figure.  They
are in $M$, but it is not at all clear if they are connected\index{connected} to
the main cardiod and period two disc of $M$ that are shown in Figure \ref{FIG_BAD_MANDEL}.  If one looks closer, one sees many more such ``dots''.
In Section \ref{SEC:MANDEL_OUTSIDE_IN} we will use a smart coloring of $\C
\setminus M$ to better understand this issue.  We will then state a theorem of Douady\index{Douady} and Hubbard\index{Hubbard},
which clears up this mystery.

The Mandelbrot set\index{Mandelbrot set} was initially discovered around 1980, but the historical details are a bit controversial.
We refer the reader to Appendix G from \cite{MILNOR} for an unbiased account.
(The reader who seeks out controversy may enjoy \cite{HORGAN}.)

The corollary to the Fatou-Julia Lemma implies that $M_0 \subset M$.  In other words,
the Mandelbrot set\index{Mandelbrot set} is an ``outer approximation'' to our set $M_0$.
The reader should compare Figure \ref{FIG_MANDEL_BW} with
Figure \ref{FIG_BAD_MANDEL} for a better appreciation of much progress we've
made!

\begin{exercise}
Prove that $M_0 \neq M$ by exhibiting a parameter\index{parameter} $c$ for which $p_c^{\circ n}(0)$ remains bounded
but with $p_c$ having no attracting\index{periodic orbit!attracting} periodic orbit.
\end{exercise}

\begin{DENSITYconjecture}
$\overline{M_0} = M$.
\end{DENSITYconjecture}

\noindent
Although this conjecture is currently unsolved, the corresponding result for real polynomials
$p_c(x) = x^2+c$ with $x,c \in \mathbb{R}$ was proved by Lyubich\index{Lyubich} \cite{LY_DENSITY}
and by Graczyk\index{Graczyk}-{\'S}wi{\c{a}}tek\index{S@{\'S}wi{\c{a}}tek} \cite{GS}.  Both proofs use complex techniques to solve
the real problem.

We have approached the definition of $M$ from an unusual perspective, i.e.
``from the inside out''.  In the next section we will use the fixed point\index{fixed point} at
$\infty$ for $p_c$ to study $M$ again, but ``from the outside in''.  It is the more traditional
way of introducing $M$.

\section*{Lecture 3: ``Complex Dynamics from the Outside In''}  
\label{SEC:MANDEL_OUTSIDE_IN}
\setcounter{section}{3}
\setcounter{subsection}{0}

\begin{definition}
The {\em filled Julia set\index{filled Julia set}} of $p_c(z) = z^2 + c$ is
\begin{align*}
K_c:= \{z \in \C \, : \, \mbox{$p_c^{\circ n}(z)$ remains bounded for all $n \geq 0$}\}.
\end{align*}
\end{definition}
If $p_c$ has an attracting\index{periodic orbit!attracting} periodic orbit
$\mathcal{O}$, then the basin of attraction $\A(\mathcal{O})$ is contained in
$K_c$.    However, $K_c$ is defined for any $c \in \C$, even if $p_c$ has no
attracting\index{periodic orbit!attracting} periodic orbit in $\C$.  

There is a natural way to extend $p_c$ as a function 
\begin{align*}
p_c : \C \cup \{\infty\} \rightarrow \C \cup \{\infty\}.
\end{align*}
 (More formally, the space $\C \cup \{\infty\}$ is called the Riemann Sphere; see \cite[Section 1.7]{SS}.)  This extension satisfies $p_c(\infty) = \infty$ and, 
by your solution to Exercise \ref{EX_ESCAPE}, $\infty$ always has a non-empty basin
of attraction:
\begin{align*}
\A(\infty) := \{z \in \C \, : \, p_c^{\circ n}(z) \rightarrow \infty\} = \C \setminus K_c.
\end{align*}
Thus, $\infty$ is an attracting fixed point of $p_c$ for any parameter $c \in \C$.  
In this way, 
the definition of $K_c$ is always related to basin of attraction for an attracting  fixed point, even if
$p_c$ has no attracting periodic point in $\C$.
A detailed study of $\A(\infty)$ will help us to prove nice theorems later in this subsection.

\begin{remark} The Fatou-Julia Lemma still applies to the extended function $p_c : \C \cup \{\infty\} \rightarrow \C \cup \{\infty\}$.  If you follow through the details of how this extension is done, you find that $\infty$ is a critical point of $p_c$ for every~$c$.
\end{remark}

\begin{definition}
The {\em Julia set\index{Julia set}} of $p_c(z) = z^2 + c$ is $J_c := \partial K_c$, the \index{boundary} boundary of $K_c$.
\end{definition}

\begin{exercise}\label{EX_INVARIANCE}
Check that for any $c \in \C$ the sets $K_c$ and $J_c$ are {\em totally invariant} meaning that
$z \in K_c \Leftrightarrow p_c(z) \in K_c$ and $z \in J_c \Leftrightarrow p_c(z) \in J_c$.
\end{exercise}

\begin{exercise}
Use the Cauchy Estimates\index{Cauchy Estimates} and invariance of $J_c$ to prove that any repelling periodic point for $p_c$ is in $J_c$.
\end{exercise}

Before drawing some computer images of Julia sets\index{Julia set}, it will be helpful to study $\A(\infty)$ a bit more.

\begin{definition}
A {\em harmonic function\index{harmonic function}} $h: \C \rightarrow \mathbb{R}$ is
a function with continuous\index{continuous} second partial derivatives $h(x+iy) \equiv h(x,y)$ satisfying
$$\frac{\partial^2 h}{\partial x^2} + \frac{\partial^2 h}{\partial y^2} = 0.$$
\end{definition}

One can use the Cauchy-Riemann Equations\index{Cauchy-Riemann Equations} to verify that the real or imaginary part of an analytic\index{analytic} function is harmonic and also that any harmonic function\index{harmonic function} can be written (locally) as the real or imaginary part
of some analytic\index{analytic} function.  Thus, there is a close parallel between the theory of analytic\index{analytic} functions
and of harmonic functions\index{harmonic function}.  We will only need two facts which follow directly from their analytic\index{analytic} counterparts:

\begin{MAXtheorem}Suppose $h(z)$ is harmonic in a domain\index{domain} $U$ and $h(z)$ achieves its maximum or minimum at a point $z_0 \in U$.
Then $h(z)$ is constant on $U$.

If, moreover, $\overline{U}$ is compact\index{compact} and $h$ extends continuously to
$\overline{U}$, then $h$ achieves its maximum and minimum on the \index{boundary} boundary of $U$.
\end{MAXtheorem}

\begin{ULHtheorem} Suppose $h_k: U \rightarrow \mathbb{R}$ is a sequence of harmonic functions\index{harmonic function}
and $h: U \rightarrow \mathbb{R}$ is some other function.
If for any compact\index{compact} $K \subset U$ we have that $\{h_k\}$ converges uniformly to $h$ on $K$, then $h$ is harmonic on $U$.

Moreover, any (repeated) partial derivative of $h_k$ converges uniformly to the corresponding partial derivative of $h$ on any compact\index{compact} $K \subset U$.
\end{ULHtheorem}

\begin{lemma} \label{LEM_GREEN}
The following limit exists
\begin{align*}G_c(z) := \lim_{n \rightarrow \infty} \frac{1}{2^n} \log_+ \left|p_c^{\circ n}(z)   \right| \quad \mbox{where} \quad \log_+(x) = {\rm max}(\log(x),0)
\end{align*}
for any parameter\index{parameter} $c \in C$ and any $z \in \C$.  For each $c$ the 
resulting function $G_c: \C \rightarrow \mathbb{R}$ is called the {\em Green function\index{Green function}} associated to $p_c$.
It satisfies:
\begin{itemize}
\item[(i)] $G_c$ is continuous\index{continuous} on $\C$ and harmonic on $\A(\infty)$,
\item[(ii)] $G_c(p_c(z)) = 2 G_c(z)$,
\item[(iii)] $G(z) \approx \log|z|$ for $|z|$ sufficiently large, and
\item[(iv)] $G(z) = 0$ iff $z \in K_c$.
\end{itemize}
\end{lemma}
\noindent
The Green Function $G_c$ is interpreted as describing the rate at which the orbit of initial condition $z_0 = z$ escapes
to infinity under iteration of $p_c$.
(The proof of Lemma \ref{LEM_GREEN} is quite similar to the proofs of K\oe nig's Theorem and B\"ottcher's Theorem, so we will omit it.)

\begin{figure}[h!]
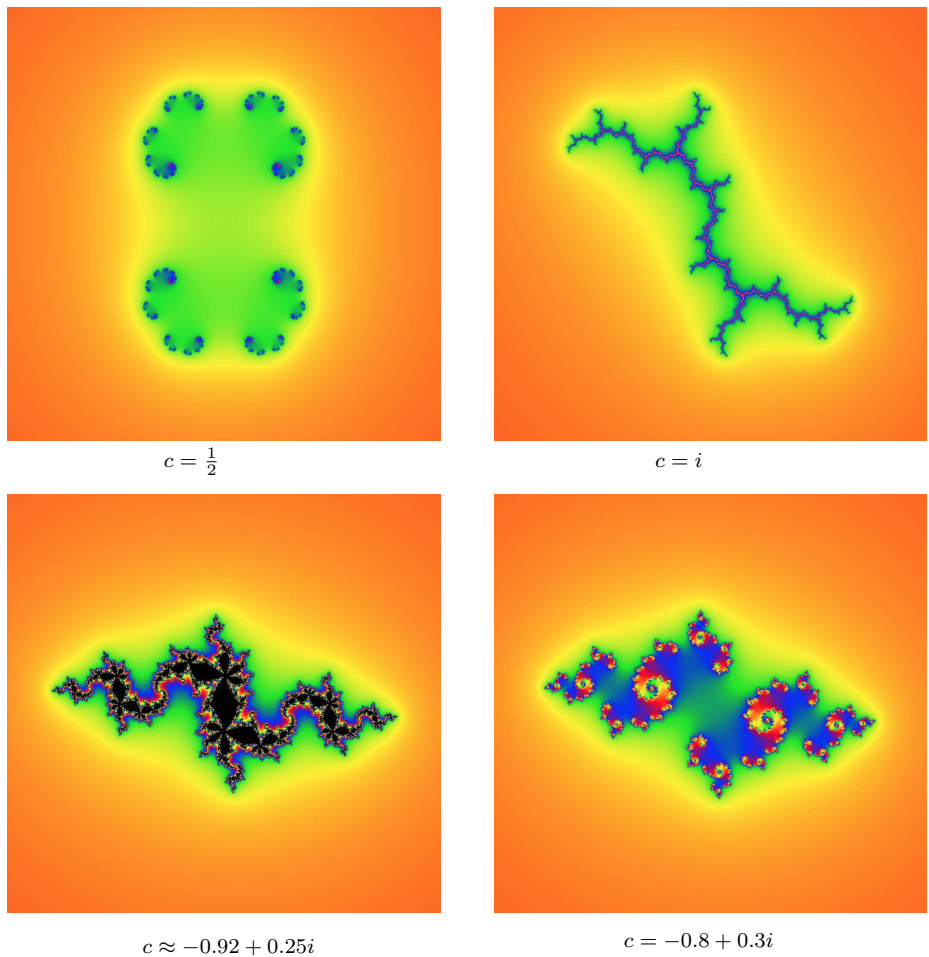
\caption{Filled Julia sets for four values of $c$.  The basin of attraction for $\infty$ is colored
according to the value of $G_c(z)$.  \label{FIG_FILLED_JULIA}}
\end{figure}

It is customary when drawing filled Julia sets\index{filled Julia set} on the computer to color
$\A(\infty) = \C \setminus K_c$ according to the values of $G_c(z)$.  This is
especially helpful for parameters\index{parameter} $c$ at which $p_c$ has no attracting\index{periodic orbit!attracting} periodic
orbit.  Using how the colors cycle one can ``view'' where $K_{c}$ should be.
In Figure \ref{FIG_FILLED_JULIA} we show the filled Julia sets\index{filled Julia set} for four
different values of $c$.  (Among these is $c=\frac{1}{2}$, from Example \ref{EXAMPLE_C_ONE_HALF}.  We can now see where
the bounded orbits are.)

For the parameter\index{parameter} values $c=\frac{i}{4}, -1,$ and $c \approx -0.12+0.75i$, the
filled Julia set\index{filled Julia set} is the closure of the basin of attracting\index{periodic orbit!attracting} periodic orbit.
Thus, \hbox{Figures \ref{FIG_BASINS1}-\ref{FIG_BASINS3}} also depict the filled
Julia sets for these parameter\index{parameter} values.

\begin{remark} Like the ancient people who named the constellations, people doing complex dynamics also have active imaginations.  They have named the filled Julia sets\index{filled Julia set} for $c = -1$ the ``basilica'' and the filled Julia set\index{filled Julia set} for $c \approx -0.12 + 0.75i$ ``Douady's\index{Douady} Rabbit''.
\end{remark}

\begin{figure}[h!]
\scalebox{1.2}{
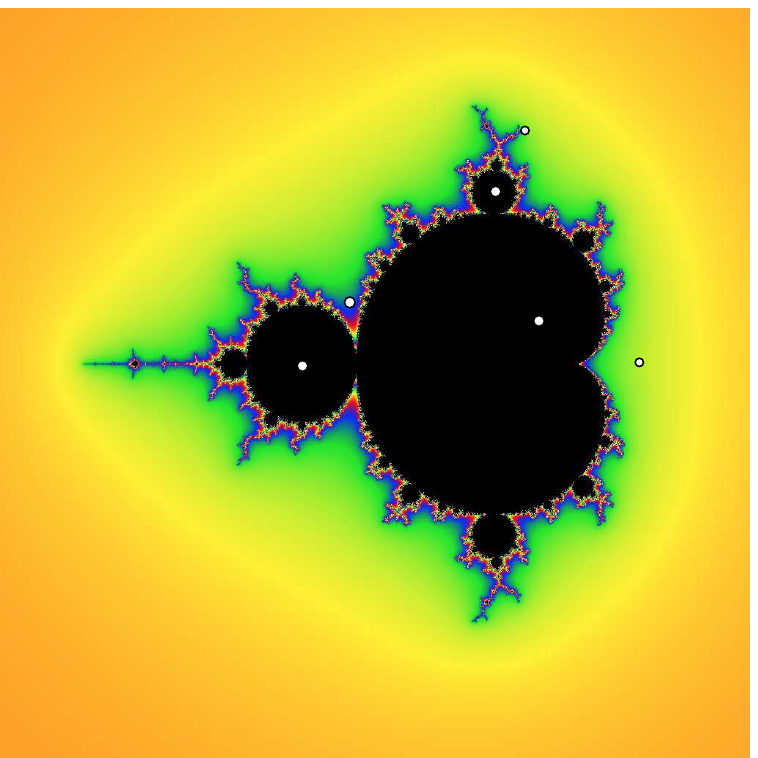
}
\caption{Mandelbrot set\index{Mandelbrot set} with the approximate locations of parameters\index{parameter} \hbox{$c = \frac{i}{4}, -1, -0.12+0.75i, i, -0.8+0.3i$,} and $\frac{1}{2}$ labeled. \label{FIG_MANDELBROT_LABELLED}}
\end{figure}

The Green function\index{Green function} also helps us to make better computer pictures of the
Mandelbrot set\index{Mandelbrot set}.  The value $G_c(0)$ expresses the rate at which the critical point\index{critical point} $0$ of $p_c$ escapes to $\infty$ under iteration of $p_c$.  Thus, points
$c$ with larger values of $G_c(0)$ should be farther away from $M$.  Therefore, it is
customary to color $\C \setminus M$ according to the values of $G_c(0)$, as in
Figure \ref{FIG_MANDELBROT_LABELLED}.  It is interesting to compare Figures \ref{FIG_MANDELBROT_LABELLED} 
and \ref{FIG_MANDEL_BW}.  It now looks more plausible that the black ``dots'' in Figure \ref{FIG_MANDEL_BW}
might be connected\index{connected} to the ``main part'' of $M$.

The Green function\index{Green function} is not only useful for making nice pictures.  It also plays a key role in the proof of:

\begin{TOP_CHAR_Mtheorem}\label{THM_MANDEL_VS_CONNECTED}
$K_c$ is connected\index{connected} if and only if $c~\in~M$.
\end{TOP_CHAR_Mtheorem}

\begin{figure}[h!]
\scalebox{1.2}{
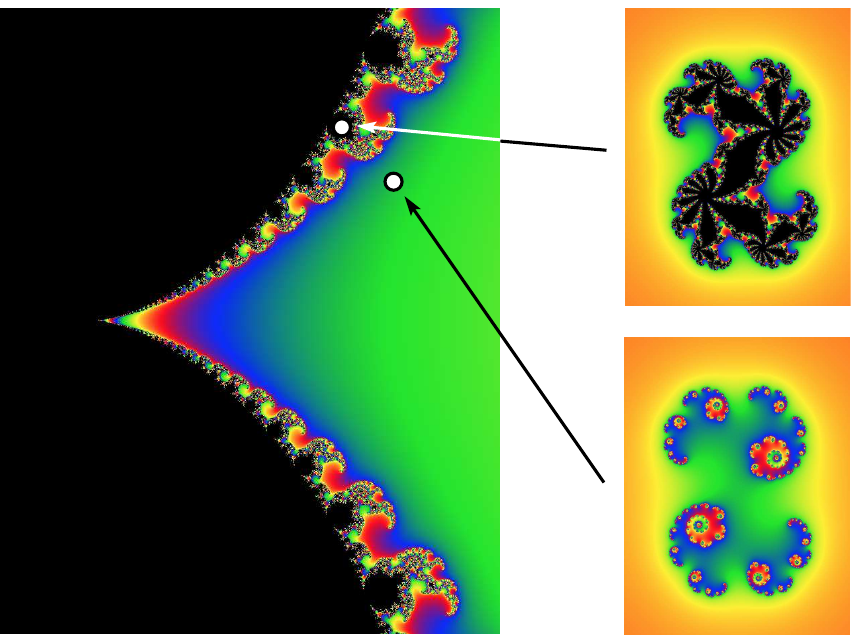
}
\caption{Left: Zoomed-in view of the Mandelbrot set\index{Mandelbrot set} near the cusp at $c=\frac{1}{4}$.  Right:
two filled Julia sets\index{filled Julia set} corresponding to points inside $M$ and outside of $M$. \label{FIG_MANDELBROT_V_JULIA}}
\end{figure}

\noindent
We illustrate this theorem with Figure
\ref{FIG_MANDELBROT_V_JULIA}.   The reader may also enjoy comparing the parameter\index{parameter} values shown in Figure \ref{FIG_MANDELBROT_LABELLED} with their Filled Julia Sets shown in previous figures.

According to the definition (\ref{EQN_DEF_MANDELBROT})  of $M$, this statement is equivalent to

\begin{TOP_CHAR_Mtheorem_Prime}\label{THM_MANDEL_VS_CONNECTED_PRIME}
$K_c$ is connected\index{connected} if and only if the orbit $\{p_c^{\circ n}(0)\}$ of the critical point\index{critical point} $0$ of $p_c$ remains bounded.
\end{TOP_CHAR_Mtheorem_Prime}

\noindent
Although the Mandelbrot set\index{Mandelbrot set} was not defined at the time of Fatou and
Julia's work (they lived from 1878-1929 and 1893-1978, respectively), the proof of the {\color{white}'}Topological Characterization of the Mandelbrot Set' is due to them.

\vspace{0.1in}

{\em Sketch of the proof:}
We will consistently identify $\C$ with $\mathbb{R}^2$ when taking partial derivatives and gradients of $G_c: \C \rightarrow \mathbb{R}$.
We claim that $G_c(z)$ has a critical point\index{critical point} at $z_0 \in \A(\infty)$ if and only
if $p^{\circ n}(z_0) = 0$ for some $n \geq 0$.  Consider the finite
approximates
\begin{align*}
G_{c,n}(z) := \frac{1}{2^n} \log_+|p_c^{\circ n}(z)|,
\end{align*}
which one can check converge uniformly to $G_c(z)$ on any compact\index{compact} subset of~$\C$.  For points $z \in \A(\infty)$ we can drop the subscript $+$ and use that
$\log|z|$ is differentiable\index{differentiable} on $\C \setminus \{0\}$.  Combined with the chain
rule, we see that $\frac{\partial}{\partial x}G_{c,n}(z) =
\frac{\partial}{\partial y}G_{c,n}(z) = 0$ if and only if $(p^{\circ n})'(z) =
0$. This holds if and only if $p^{\circ m}(z) = 0$ for some $0 \leq m \leq
n-1$.  The claim then follows from the Uniform Limits of Harmonic Functions\index{Uniform Limits of Harmonic Functions}
Theorem.

Suppose that the critical point\index{critical point} $0$ has bounded orbit under $p_c$.  Then,
according to the previous paragraph, $G_c$ has no critical points\index{critical point} in
$\A(\infty)$.
For any $t > 0$ let 
\begin{align*}
L_t := \{z \in \C \, : \, G_c(z) \leq t\}.
\end{align*}
 By definition, if $t < s$ then $L_t \subset L_s$.
For each $t > 0$, $L_t$ is closed and bounded since $G_c: \C \rightarrow \mathbb{R}$ is
continuous\index{continuous} and $G_c(z) \rightarrow \infty$ as $|z| \rightarrow \infty$, respectively.  Therefore, by the
Heine-Borel Theorem\index{Heine-Borel Theorem}, $L_t$ is compact\index{compact}.  Since $K_c \neq \emptyset$ and $G_c(z) = 0$ on $K_c$, each $L_t$ is non-empty.  

Since $z \in K_c$ if and only if $G_c(z) = 0$, we can write $K_c$ as a nested intersection of non-empty compact\index{compact} sets:
\begin{align*}
K_c = \bigcap_{n \geq 1} L_{1/n}.
\end{align*}
If we can show that $L_t$ is connected\index{connected} for each $t > 0$, then Exercise \ref{EX_NESTED} will imply that
$K_c$ is connected\index{connected}. 

Since $G_c(z) \approx \log|z|$ for $|z|$ sufficiently large, there exists $t_0
> 0$ sufficiently large so that $L_{t_0}$ is connected\index{connected} (it is almost a
closed disc\index{closed disc} of radius $\log t_0$).  We will show that for any $0 < t_1 < t_0$
the sets $L_{t_1}$ and $L_{t_0}$ are homeomorphic (i.e, there is a continuous\index{continuous}
bijection with continuous\index{continuous} inverse from $L_{t_1}$ to $L_{t_0}$).  Since $L_{t_0}$ is connected\index{connected},
this will imply that $L_{t_1}$ is also connected\index{connected}.

The following is a standard construction from Morse Theory; see \cite[Theorem 3.1]{MIL_MORSE}.
Because $G_c$ is harmonic and has no critical points\index{critical point}  on $\A(\infty)$,
 $-\nabla G_c$
is a non-vanishing smooth vector field on $\A(\infty)$.
It is a relatively standard smoothing construction to define a new vector field $\bm{V}: 
\mathbb{R}^2 \rightarrow \mathbb{R}^2$ that is smooth on all of $\C \equiv
\mathbb{R}^2$ and equals $\frac{-\nabla G_c}{\|\nabla G_c\|^2}$ for $z \in \C \setminus L_{t_1/2}$.

For any $t \in [0,\infty)$ let $\Phi_t: \mathbb{R}^2 \rightarrow \mathbb{R}^2$ denote the flow obtained by integrating
$\bm{V}$.  According to the existence and uniqueness theorem for ordinary differential equations (see, e.g., \cite{PERKO}),
$\Phi_t: \mathbb{R}^2 \rightarrow \mathbb{R}^2$ is a homeomorphism for each $t \in [0,\infty)$.  (We're using that $\bm{V}$ ``points inward'' from $\infty$ so that the solutions exist for all time.)

For any $z_0 \in \C \setminus L_{t_1/2}$ and any $0 \leq t \leq G_c(z_0) - t_1 /2$  we have
\begin{align*}
\frac{d}{dt} G_c(\Phi_t(z_0)) &= \nabla G_c(\Phi_t(z_0)) \cdot
\frac{d}{dt}\Phi_t(z_0)  =  \nabla G_c(\Phi_t(z_0)) \cdot {\bm
V}(\Phi_t(z_0))\\ &= \nabla G_c(\Phi_t(z_0)) \cdot  \frac{-\nabla
G_c(\Phi_t(z_0))}{\|\nabla G_c (\Phi_t(z_0))\|^2} = -1.
\end{align*}
In particular, 
$\Phi_{t_0-t_1} (L_{t_0}) = L_{t_1}$, implying that $L_{t_0}$ is homeomorphic to $L_{t_1}$.

\vspace{0.1in}

Now suppose that $0$ has unbounded orbit under $p_c$.  In this case, $0$ and all of its iterated preimages
under $p_c$ are \index{critical point} critical points of $G_c$.
Since $p_c$ has a simple critical point\index{critical point} at $0$, one can check
that these critical points\index{critical point} of $G_c$ are all ``simple'' in that the Hessian
matrix of second derivatives has non-zero determinant.  Moreover, by the
Maximum Principle\index{Maximum Principle}, they cannot be local minima or local maxima.  They are
therefore saddle points.  From the property $G_c(p(z)) = 2 G_c(z)$, the saddle
point at $z = 0$ is the one with the largest value of $G_c$.

There are two paths along which we can start at $0$ and walk uphill in the
steepest way possible---call them $\gamma_1$ and $\gamma_2$.  Since $0$ is the
highest critical point\index{critical point}, they lead all the way from $0$ out to $\infty$.
Together with $0$, these two paths divide $\C$ into two domains\index{domain} $U_1$ and
$U_2$.  Meanwhile, there are two directions that one can walk downhill from a
saddle point.  Walking the fastest way downhill leads to two paths $\eta_1$ and
$\eta_2$ which lead to points in $U_1$ and in $U_2$ along which $G_c(z) <
G_c(0)$.  

To make this idea rigorous, one considers the flow associated to the vector
field $-\nabla G_c$.  The saddle point $0$ becomes a saddle type fixed point\index{fixed point}
for the flow with the paths $\gamma_1$ and $\gamma_2$ being the stable manifold of
this fixed point\index{fixed point}.  The paths $\eta_1$ and $\eta_2$ are the unstable manifolds
of this fixed point\index{fixed point}. (See again~\cite{PERKO}.)

The union $\gamma_1 \cup \gamma_2 \cup \{0\}$ divides
the complex plane into two domains\index{domain} $U_1$ and $U_2$ with $\eta_1 \subset U_1$
and $\eta_2 \subset U_2$.  We claim that both of these domains\index{domain} contain points
of $K_c$.  Suppose for contradiction that one of them (say $U_1)$ does not.
Then, $U_1 \subset \A(\infty)$ and hence $G_c$ would be harmonic on $U_1$.
However, $G_c(z) \sim \log|z|$ for $|z|$ large and $G_c(z) > G_c(0)$ for points
$z \in \gamma_1 \cup \gamma_2$.  Since $G_c(z) < G_c(0)$ for points on
$\eta_1$, this would violate the Maximum Principle\index{Maximum Principle}.

\qed

\begin{remark}
A stronger statement actually holds: if $K_c$ is
disconnected\index{disconnected}, then it is a Cantor Set. (See \cite{EDGAR}
for the definition of Cantor Set.)  In particular, it is {\em totally
disconnected}: for any $z, w \in K_c$ there exist open sets $U, V \subset \C$
such that $K_c \subset U \cup V$, $z \in U$, $w \in V$, and $U \cap V =
\emptyset$.  This follows from the fact that once $G_c$ has the critical point $0 \in \A(\infty)$ then
it actually has infinitely many critical points in $\A(\infty)$.  These additional critical
points of $G_c$ are the iterated preimages of $0$ under $p_c$.

For a somewhat different proof from the one presented above, including a proof of this stronger statement, see \cite{DEV1,DEV2}.
\end{remark}

\begin{exercise} Prove that if $c \neq 0$ then $\log|z^2+c|$ has a saddle-type critical point\index{critical point} at $z=0$.

\noindent
{\em Hint:} Write $z=x+iy$ and $c=a+bi$ and use that $$\log|z^2+c| = \frac{1}{2} \log\Big(z^2+c\Big)\Big(\overline{z^2+c}\Big).$$
\end{exercise}

We will now state (without proofs) several interesting properties of the Mandelbrot set\index{Mandelbrot set}:

\begin{theorem*}{\bf (Douady\index{Douady}-Hubbard\index{Hubbard} \cite{DH})} The Mandelbrot set\index{Mandelbrot set} $M$ is connected\index{connected}.
\end{theorem*}
\noindent
(Nessim Sibony gave an alternate proof around the same time.)
This theorem clears up the mystery about the black ``dots'' in Figure \ref{FIG_MANDEL_BW}.

The following very challenging extended exercise leads the reader through a proof that $M$ is connected that is related
to the coloring of $\C \setminus M$ according to the value of $G_c(0)$.  
(It will be somewhat more convenient to consider $G_c(c) = G_c(p_c(0)) = 2G_c(0)$.)

\begin{extendedexercise}
Let $H: \C \rightarrow \mathbb{R}$ be given by $H(c) = G_c(c)$. Prove that
\begin{enumerate}
\item $H$ is continuous, 
\item $H$ is harmonic on $\C \setminus M$,
\item $H$ is identically $0$ on $M$, 
\item $\lim_{|c| \rightarrow \infty} H(c) = \infty$, and
\item  $H$ has no critical points in $\C \setminus M$. 
\end{enumerate}
(Step 5 is the hardest part.)
Use these facts to adapt the proof of the topological characterization of the Mandelbrot Set to prove that $M$ is connected.
\end{extendedexercise}

\begin{figure}[h!]
\scalebox{0.55}{
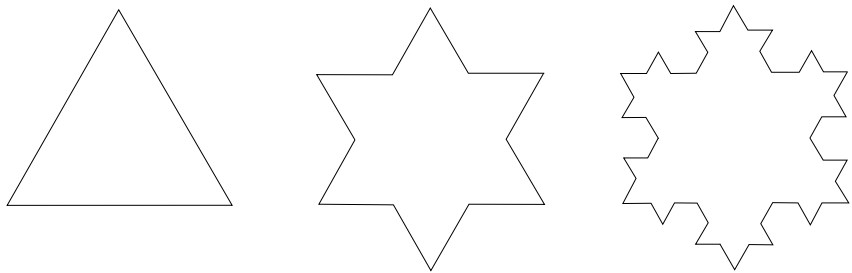
}
\caption{The Koch Curve.
\label{FIG_KOCH_CURVE}}
\end{figure}

{\em Hausdorff Dimension}\index{Hausdorff Dimension} extends the
classical notion of dimension from lines and planes to more general metric
spaces.  As the formal definition is a bit complicated, we instead illustrate the
notion with a few examples.  A line has Hausdorff dimension equal to $1$ and
the plane has Hausdorff dimension equal to $2$.    A contour \index{contour}
has Hausdorff dimension equal to $1$ because, if you zoom
in sufficiently far near any of the smooth points,
the contour appears more and more like a straight line.  However,
sets of a ``fractal nature'' can have non-integer Hausdorff Dimension.
One example is the Koch Curve, which a simple closed curve in the
plane that is obtained as the limit of the iterative process shown in Figure
\ref{FIG_KOCH_CURVE}.  No matter how far you zoom in, the Koch Curve 
looks the
same as a larger copy of itself, and not like a line!  This results in the Koch Curve having
Hausdorff dimension equal to $\log(4)/\log(3) \approx 1.26$.
We refer the reader to \cite{EDGAR} for a gentle introduction to Hausdorff Dimension.

If $S \subset \C$ contains an open subset of $\C$, then it is easy to see
that it has Hausdorff Dimension equal to $2$.  It is much harder to imagine
a subset of $\C$ that contains no such open set having Hausdorff Dimension $2$.  Therefore,  
the following theorem shows that the \index{boundary} boundary $\partial M$ of the Mandelbrot
set\index{Mandelbrot set} $M$ has amazing complexity.  It also shows that for
many parameters\index{parameter} $c$ from $\partial M$ the Julia
set\index{Julia set} $J_c$ has amazing complexity.

\begin{theorem*}{\bf (Shishikura\index{Shishikura} \cite{SHISHIKURA})} The \index{boundary} boundary of the
Mandelbrot set\index{Mandelbrot set} $\partial M$ has Hausdorff dimension equal to $2$.  Moreover,
for a dense set of parameters\index{parameter} $c$ from the boundary of $M$ the Julia set\index{Julia set} $J_c$
has Hausdorff Dimension equal to~$2$.  
\end{theorem*}

\begin{figure}[h!]
\scalebox{1.4}{
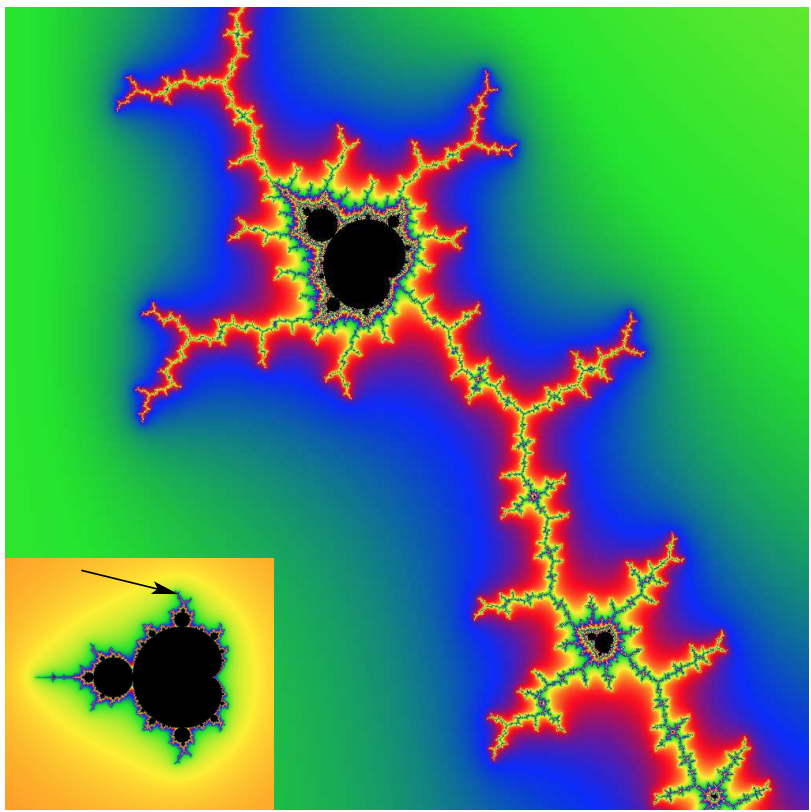
}
\caption{Zoomed-in view of part of the Mandelbrot set\index{Mandelbrot set}  showing two smaller copies.
The approximate location where we have zoomed in is marked by the tip of the arrow in the inset figure.
\label{FIG_SMALL_MANDELBROT}}
\end{figure}

Another interesting property of the Mandelbrot set\index{Mandelbrot set} is the appearance of ``small
copies'' within itself.  (Some of these were the ``dots'' from Figure~\ref{FIG_MANDEL_BW}.)  Figure \ref{FIG_SMALL_MANDELBROT} shows a zoomed in view
of $M$, where several small copies of $M$ are visible.  These copies are
explained by the renormalization theory \cite{DH_RENORM,LY_BABY_MANDELBROT}.

It would be remiss to not include one of the most famous conjectures about the Mandelbrot set\index{Mandelbrot set}.
We first need
\begin{definition}
A topological space $X$ is {\em locally connected\index{connected!locally connected}} if for every point $x \in X$
and any open set $V \subset X$ that contains $x$ there is another connected\index{connected} open set $U$ with
$x \in U \subset V$.
\end{definition}

\begin{MLCconjecture}
The Mandelbrot set\index{Mandelbrot set} $M$ is locally connected\index{connected!locally connected}.
\end{MLCconjecture}

\noindent
According to the Orsay Notes\index{Orsay Notes} \cite{ORSAY} of Douady\index{Douady} and Hubbard\index{Hubbard}, if this were the case, then one could
have a very nice combinatorial description of $M$.  Given a proposed way that $p_c$ acts 
on the Julia set\index{Julia set} $J_c$ (described by means of the so called Hubbard Tree) one
can use this combinatorial description of $M$ to find the desired value of~$c$.

\begin{figure}[h!]
\scalebox{0.9}{
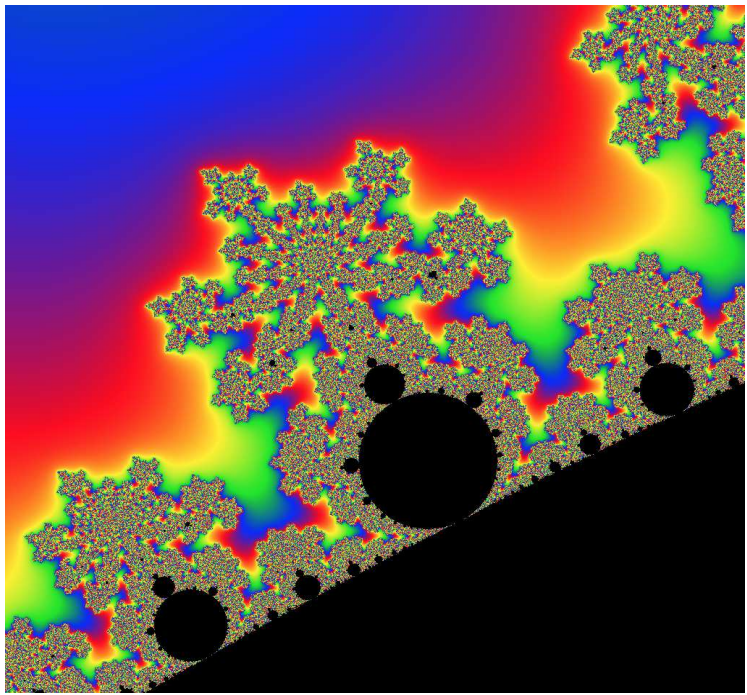
}
\caption{Another zoomed-in view of part of the Mandelbrot set\index{Mandelbrot set}.
\label{FIG_MANDEL_ZOOMED}}
\end{figure}

To better appreciate the difficulty in proving 
the MLC Conjecture, we include one more zoomed-in image of the Mandelbrot set\index{Mandelbrot set} in Figure \ref{FIG_MANDEL_ZOOMED}.

\vspace{0.1in}

Let us finish the section, and our discussion of iterating quadratic
polynomials, by returning to mathematics that can be done by undergraduates.
The reader is now ready to answer Question 4 from Section \ref{SEC_QUESTIONS}:

\begin{extendedexercise}
Prove that  for every $m \geq 1$ there exists a parameter\index{parameter} $c \in \C$ such that $p_c(z)$ has an attracting\index{periodic orbit!attracting} periodic
orbit of period exactly~$m$.

\noindent
{\em Hint:} prove that there is a parameter\index{parameter} $c$ such that $p_c^{\circ m}(0) = 0$ and
$p_c^{\circ j}(0) \neq 0$ for each $0 \leq j < m$.
\end{extendedexercise}

\section*{Lecture 4: ``Complex dynamics and astrophysics.''}
\setcounter{section}{4}
\setcounter{subsection}{0}

Most of the results discussed in Sections 1-3 of these notes are now quite
classical.   Let us finish our lectures with a beautiful and quite modern application
of the Fatou-Julia Lemma\index{Fatou-Julia Lemma} to a problem in astrophysics \cite{KS,KN}.  We also
mention that there are connections between complex dynamics and the Ising model from statistical
physics (see \cite{BLR1,BLR2} and the references therein) and the study of
droplets in a Coulomb gas \cite{LM1,LM2}.

\subsection{Gravitational Lensing}

\begin{figure}[h!]
\begin{center}
\scalebox{0.8}{
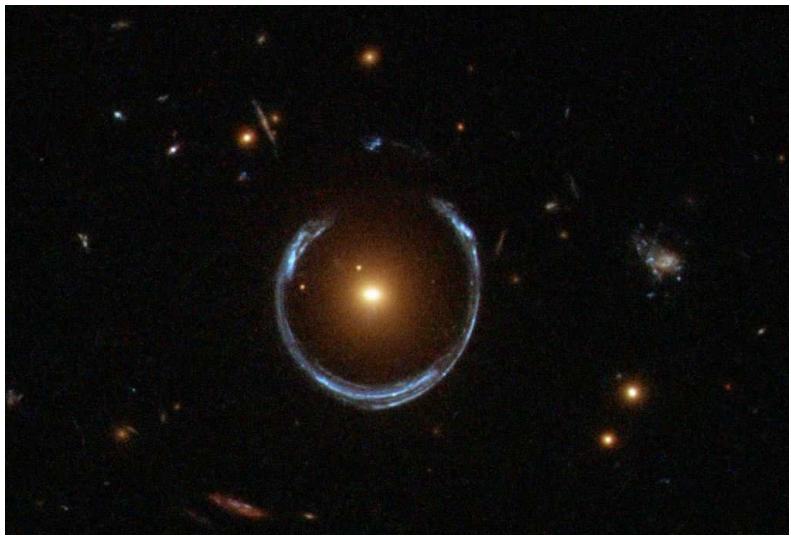%
}
\caption{An Einstein Ring.  For more information, see \protect\url{http://apod.nasa.gov/apod/ap111221.html}.
 \label{FIG_NASA_IMAGE1}}
\end{center}
\end{figure}

\begin{figure}[h!]
\begin{center}
\scalebox{0.8}{
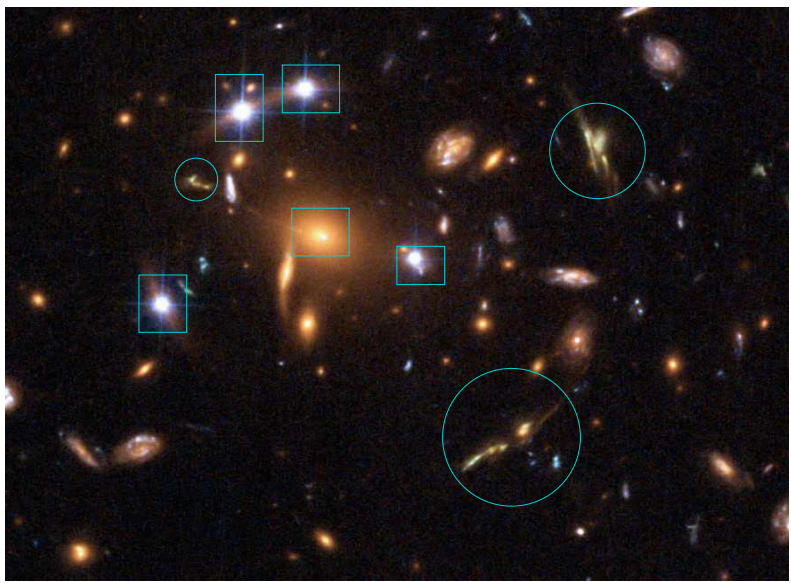%
}
\caption{
Five images of the same quasar
(boxed) and three images of the same galaxy (circled).  The middle image of the quasar (boxed) is behind the small galaxy that does the lensing.  For more information, see  \protect\url{http://www.nasa.gov/multimedia/imagegallery/image_feature_575.html}
 \label{FIG_NASA_IMAGE2}}
\end{center}
\end{figure}

Einstein's Theory of General Relativity predicts that if a point mass is placed directly between
an observer and a light source, then the observer will see a ring of light, called an ``Einstein Ring''.
The Hubble Space Telescope has sufficient power to see these rings---one such image is 
shown in Figure \ref{FIG_NASA_IMAGE1}.
  If the point mass
is moved slightly, the observer will see two different images of the same light source.  With more complicated distributions
of mass, like $n$ point masses, the observer can see more complicated images, resulting from a single point light source.
Such an image is shown in Figure \ref{FIG_NASA_IMAGE2}.  (Thanks to NASA for these images and their interpretations.)

There are many excellent surveys on gravitational lensing that are written for the
mathematically inclined reader, including \cite{KN_SURVEY,PETTERS_SURVEY,STRAUMANN}, as
well as the book~\cite{PETTERS_BOOK}.  We will be far more brief, with the goal of this lecture being
to explain how Rhie\index{Rhie}~\cite{RHIE} and Khavinson\index{Khavinson}-Neumann\index{Neumann} \cite{KN} answered the question:\\

\noindent
{\it 
What is the maximum number of images that a single light source can have when lensed by $n$ point masses? }\\

\noindent
We will tell some of the history of how this problem was solved and then focus on the role played by the Fatou-Julia Lemma\index{Fatou-Julia Lemma}.

\vspace{0.1in}

Suppose that $n$ point masses lie on a plane that is nearly perpendicular to
the line of sight between the observer and the light source and that they lie relatively close to the line of sight.
If we describe their positions relative to the line of sight to the light source by complex numbers $z_j$ and their normalized masses by $\sigma_j > 0$ for $1\leq j \leq n$,
then the images of the light source seen by the observer are given by  solutions $z$ to the {\em Lens Equation:}\index{Lens Equation}
\begin{equation}\label{LENS EQUATION}
z = \sum_{j=1}^{n} \frac{\sigma_j}{\bar{z} - \bar{z_j}}.
\end{equation}
The ``mysterious'' appearance of complex conjugates on the right hand side of this equation makes
it difficult to study.
It will be explained in Section~\ref{SUBSEC_DERIVE_LENS_EQ}, where we derive (\ref{LENS EQUATION}).

\begin{exercise} Verify that (\ref{LENS EQUATION}) gives a full circle of solutions (Einstein Ring) when there is just one mass at $z_1 = 0$.  Then, verify that when $z_1 \neq 0$ there are two solutions.
Can you find a configuration of two masses so that (\ref{LENS EQUATION}) has five solutions?
\end{exercise}

\begin{remark} Techniques from complex analysis extend nicely to lensing by mass distributions more complicated than
finitely many points, including elliptical \cite{FKK} and spiral \cite{BEFKKL} galaxies.    
\end{remark}

The right hand side of (\ref{LENS EQUATION}) is of the form $\bar{r(z)}$, where $r(z)$ is a rational function $r(z) = \frac{p(z)}{q(z)}$ of degree $n$.  (The degree of a rational function is the maximum of the degrees of its numerator and denominator.)
Thus, our physical question becomes the problem of bounding the number of solutions 
to an equation of the form 
\begin{align}\label{EQN_RATIONAL_WEIRD}
z = \bar{r(z)} 
\end{align}
in terms of $n = \deg(r(z))$.
Sadly, the Fundamental Theorem of Algebra\index{Fundamental Theorem of Algebra} cannot be applied to
\begin{align}\label{EQN_REDUCED}
z \bar{q(z)} - \bar{p(z)} = 0
\end{align}
because the resulting equation is a polynomial in {\em both} $z$ and $\bar{z}$.  If one writes $z=x+iy$ with $x,y \in \mathbb{R}$, one can 
change (\ref{EQN_REDUCED}) to a system of two real polynomial equations
\begin{align*}
a(x,y) &:= \Re\left(z \, \bar{q(z)} - \bar{p(z)} \right) = 0 \quad \mbox{and} \\ 
b(x,y) &:=\Im\left(z \, \bar{q(z)} - \bar{p(z)}\right) = 0,
\end{align*}
each of which has degree $n+1$.  So long as there are no curves of common zeros
for $a(x,y)$ and $b(x,y)$, Bezout's Theorem\index{Bezout's Theorem} (see, e.g., \cite{KIRWAN}) gives a
bound on the number of solutions by $(n+1)^2$.  

In 1997, Mao, Petters, and Witt \cite{MPW} exhibited configurations of $n$
point masses at the vertices of a regular polygon in such a way that $3n+1$
solutions were found.  They conjectured a linear bound for the number of
solutions to~(\ref{LENS EQUATION}).
For large $n$ this
would be significantly better than the bound given by Bezout's Theorem.

In 2003, Rhie\index{Rhie} \cite{RHIE} showed that if one takes the configuration of masses
considered by Mao, Petters, and Witt and places a sufficiently small mass
centered at the origin, then one finds $5n-5$ solutions to (\ref{LENS
EQUATION}).  (We refer the reader also to \cite[Section 5]{BHJR} for an
another exposition on Rhie\index{Rhie}'s examples.)

In order to address a problem on harmonic mappings\index{harmonic mapping} $\C \rightarrow \C$ posed by Wilmshurst in
\cite{WILMS}, 
in 2003 Khavinson\index{Khavinson} and {\'S}wi{\c{a}}tek\index{S@{\'S}wi{\c{a}}tek} studied the number of solutions to 
$z = \bar{p(z)}$ where $p(z)$ is a complex polynomial\index{complex polynomial}.  They proved

\begin{theorem*}{\bf (Khavinson-{\'S}wi{\c{a}}tek \cite{KS})}\label{THM_KS}
Let $p(z)$ be a complex polynomial\index{complex polynomial} of degree $n \geq 2$.  Then, $z=\overline{p(z)}$ has at most $3n-2$ solutions.
\end{theorem*}

Khavinson\index{Khavinson} and Neumann\index{Neumann} later adapted the techniques from \cite{KS} to prove

\begin{theorem*}{\bf (Khavinson-Neumann \cite{KN})}
Let $r(z)$ be a rational function of degree $n \geq 2$.  Then, $z = \overline{r(z)}$ has at most $5n-5$ solutions.
\end{theorem*}

\noindent
Apparently, Khavinson\index{Khavinson} and Neumann\index{Neumann} solved this problem because of its mathematical interest and only later were informed
that they had actually completed the solution to our main question of this lecture:\\

\noindent
{\it
When lensed by $n$ point masses, a single light source can have at most $5n-5$ images.}\\

\begin{remark}
Geyer \cite{G} used a powerful theorem of Thurston to show that for every $n
\geq 2$ there is a polynomial $p(z)$ for which $z = \bar{p(z)}$ has $3n-2$
solutions, thus showing that Theorem \ref{THM_KS} is sharp.  It would be interesting to
see an ``elementary'' proof.
\end{remark}

\subsection{Sketching the proof of the $5n-5$ bound}

We provide a brief sketch of the proof of the Khavinson\index{Khavinson}-Neumann\index{Neumann} upper bound in the special case that
\begin{align}\label{EQN_FORMAT_FOR_R}
r(z) = \sum_{j=1}^{n} \frac{\sigma_j}{z - z_j},
\end{align}
with each $\sigma_j > 0$.
It is the case arising in the Lens Equation\index{Lens Equation} (\ref{LENS EQUATION}).
The locations $\{z_1,\ldots,z_n\}$ of the masses are called {\em poles} of $r(z)$.  They satisfy $\lim_{z \rightarrow z_j} |r(z)| = \infty$ for any $1 \leq j \leq n$.

The function 
\begin{align}
f: \C \setminus \{z_1,\ldots,z_n\} \rightarrow \C \quad \mbox{given by} \quad f(z) = z - \bar{r(z)}
\end{align}
is an example of a {\em
harmonic mapping\index{harmonic mapping} with poles} since its real and imaginary parts are harmonic.  
It is orientation preserving near a point $z_\bullet$ with $|r'(z_\bullet)| < 1$ and orientation reversing (like a reflection $z \mapsto \bar{z}$) near points with $|r'(z_\bullet)| > 1$.
A zero $z_\bullet$ of $f$ is {\em simple} if $|r'(z_\bullet)| \neq 1$ and a simple zero is called {\em sense preserving} if $|r'(z_\bullet)| < 1$ and {\em sense-reversing} if $|r'(z_\bullet)| > 1$.

\vspace{0.1in}
\noindent
{\bf Step 1: Reduction to Simple Zeros.}  Suppose $r(z)$ is of the form~(\ref{EQN_FORMAT_FOR_R}) and $f(z) = z-\bar{r(z)}$ has $k$ zeros,
some of which are not simple.  Then, one can show that there is an arbitrarily small perturbation of the locations
of the masses so that the resulting rational function $s(z)$ produces $g(z) = z - \bar{s(z)}$ having at least as many zeros as $f(z)$ all of which are simple.

Therefore, it suffices to consider rational functions
$r(z)$ of the form (\ref{EQN_FORMAT_FOR_R}) such that each zero of $f(z) = z-\bar{r(z)}$ is simple.

\vspace{0.1in}
\noindent
{\bf Step 2: Argument Principle\index{Argument Principle} for Harmonic Mappings.} Suffridge and Thompson \cite{ST}
adapted the Argument Principle\index{Argument Principle} to harmonic mappings\index{harmonic mapping} with poles $f: \C \setminus
\{z_1,\ldots,z_k\} \rightarrow \C$.  

Since $r(z)$ has the form (\ref{EQN_FORMAT_FOR_R}), $\lim_{z \rightarrow
\infty} |r(z)| = 0$.  Therefore, we can choose $R > 0$ sufficiently large so
that all of the poles of $r(z)$ lie in $D(0,R)$ and the change of argument for
$f(z)$ while traversing $\gamma = \partial D(0,R)$ counter clockwise is $1$.
This variant of the argument principle then gives 
\begin{align}\label{EQN_ARG_PRINCE}
(m_+ - m_-) + n = 1
\end{align}
where $m_+$ is the number of sense preserving zeros, $m_-$ is the number of sense reversing zeros, and $n$ is the number of poles.  (We are using that all of the zeros and poles are simple so that they do not need to be counted with multiplicities.)

\vspace{0.1in}
\noindent
{\bf Step 3: Fatou-Julia Lemma\index{Fatou-Julia Lemma} Bound on  $\bm m_+$:}
Zeros of $f(z)$ correspond to fixed points\index{fixed point} for the anti-analytic\index{anti-analytic} mapping
\begin{align*}
z \mapsto \bar{r(z)}.
\end{align*}
Moreover, sense preserving zeros correspond to attracting\index{periodic orbit!attracting} fixed points\index{fixed point} (those with $|r'(z_\bullet)| < 1$).

Since the coefficients of $r(z)$ are real, taking the second iterate yields
\begin{align*}
Q(z) = \bar{r\left(\bar{r(z)}\right)} = r(r(z)),
\end{align*}
which is an analytic\index{analytic} rational mapping of degree $n^2$.  Such a mapping has $2n^2 -2$ critical points\index{critical point} and an adaptation
of the Fatou-Julia Lemma\index{Fatou-Julia Lemma} implies that each attracting\index{periodic orbit!attracting} fixed point\index{fixed point} of $Q$ attracts a critical point\index{critical point}.  

However, the chain rule gives that critical points\index{critical point} of $Q(z) = r(r(z))$ are the
critical points\index{critical point} of $r(z)$ and their inverse images under $r(z)$.  Since a generic point has $n$ inverse images under $r$,
this can be used to show that each attracting\index{periodic orbit!attracting} fixed point\index{fixed point} of $Q(z)$ actually attracts $n+1$ critical points\index{critical point} of~$Q$.  Therefore,
$Q(z)$ has at most $2n-2$ attracting\index{periodic orbit!attracting} fixed points\index{fixed point}.

Since any sense preserving zero for $f(z)$ is an attracting\index{periodic orbit!attracting} fixed point\index{fixed point} for~$Q$, we conclude that $m_+ \leq 2n-2$.

\vspace{0.1in}
\noindent 
{\bf Step 4: Completing the proof:}
Since $m_+ \leq 2n-2$, Equation (\ref{EQN_ARG_PRINCE}) implies $m_- \leq 3n-3$.  Therefore, the total number of zeros is 
\begin{align*}
m_+ + m_- \leq 5n-5.
\end{align*}

The reader is encouraged to see \cite{KN} for the full details, including how to prove the bound for general rational functions $r(z)$.

\subsection{Derivation of the Lens Equation\index{Lens Equation}}
\label{SUBSEC_DERIVE_LENS_EQ}

This derivation is a synthesis
of ideas from \cite{ASADA} and  \cite[Section 3.1]{PETTERS} that was written jointly with Bleher,
Homma, and Ji when preparing \cite{BHJR}.  Since it was not included in the published version of \cite{BHJR}, we present it here.

\begin{figure}
\begin{center}
\scalebox{0.7}{
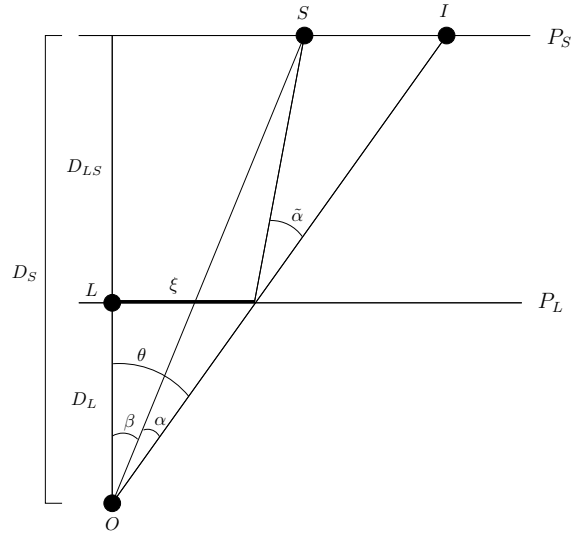%
}
\caption{$S$ is the light source, $I$ is an image, $O$ is the observer, $L$ is a point mass, $P_L$ is the lens plane, $P_S$ is the source plane.}
\label{LENS DIAGRAM}
\end{center}
\end{figure}

We will first derive the Lens Equation\index{Lens Equation} for one point mass using Figure~\ref{LENS DIAGRAM}, and then adapt it to $N$ point masses. Suppose the observer is located at point $O$, the light source at a point $S$, and a mass $M$ at point $L$. Also, suppose $P_L$ is the plane perpendicular to $OL$ that contains $L$, and $P_S$ is the plane perpendicular to $OL$ that contains $S$. Due to the point mass, an image, $I$, will be created at angle $\alpha$ with respect to~$S$.

Einstein derived using General Relativity that the bending angle is
\begin{align}\label{EQN_ALPHA_TILDE1}
	\tilde{\alpha}=\frac{4GM}{c^2\xi},
\end{align}
where $G$ is the universal gravitational constant and $c$ is the speed of light, see ~\cite{ASADA}. 

The observer $O$ describes the location $S$ of the light source using an angle
$\beta$ and the perceived location $I$ using another angle $\theta$ (see
Figure~\ref{LENS DIAGRAM}).   By a small angle approximation, $\xi = D_L \theta$, which we substitute into (\ref{EQN_ALPHA_TILDE1}) obtaining
\begin{align}\label{EQN_ALPHA_TILDE2}
        \tilde{\alpha}(\theta)=\frac{4GM}{c^2 D_L \theta}.
\end{align}
A small angle approximation also gives that $D_{SI} = D_{LS} \tilde{\alpha}(\theta) = D_S \alpha(\theta)$.
Substituting this into $\beta=\theta-\alpha(\theta)$ gives 
\begin{equation}\label{EINSTEIN RING}
	\beta=\theta-\frac{D_{LS}}{D_S D_L}\cdot\frac{4GM}{c^2\theta}.
\end{equation}

\noindent For $\beta\neq 0$, exactly two images are produced. When $\beta=0$,
the system is rotationally symmetric about $OL$, thereby producing an
Einstein Ring, whose angular radius is given by
Equation~\eqref{EINSTEIN RING}.

In order to describe systems of two or more masses, we need to describe
locations in the source plane $P_S$ and the lens plane $P_L$ using
two-dimensional vectors of angles (polar and azimuthal angles) as observed from
$O$.  
Complex numbers will be a good way to do this:
\begin{align*}
\alpha = \alpha^{(1)} + i \alpha^{(2)}, \,\,  \tilde{\alpha} = \tilde{\alpha}^{(1)}+ i\tilde{\alpha}^{(2)}, \,\,
\beta = \beta^{(1)}+i\beta^{(2)}, \,\,  \mbox{and} \,\,  \theta = \theta^{(1)} + i \theta^{(2)}.
\end{align*}
When there is only one mass, the whole configuration must still lie in one plane, as in Figure \ref{LENS DIAGRAM}.  In particular, all four complex
numbers have the same argument, forcing us to replace the $\theta$ on the right hand side of (\ref{EQN_ALPHA_TILDE2}) with  $\bar{\theta}$:
\begin{align}\label{EQN_CONJUGATE}
        \tilde{\alpha}(\theta)=\frac{4GM}{c^2 D_L \bar{\theta} } \qquad \mbox{and hence} \qquad \beta=\theta-\frac{D_{LS}}{D_S D_L}\cdot\frac{4GM}{c^2 \bar{\theta}}.
\end{align}
{\it This is why the complex conjugate arises in the Lens Equation\index{Lens Equation} (\ref{LENS EQUATION}).}

\vspace{0.1in}

We now generalize to $n$ point masses.  Let $L$ be the
center of mass of the $n$ masses, and redefine $S_L$ as the plane that is
perpendicular to $OL$ and contains $L$.  We assume that the distance between
$L$ and the individual point masses is extremely small with respect to the
pairwise distances between $O$, $P_L$, and $P_S$. Now consider the projection
of the $n$ point masses onto $S_L$.  We continue to let $\beta \in \C$ describe the location of the center
of mass and we describe
the location of the $j^{th}$ point mass by
$\epsilon_j = \epsilon_j^{(1)} + i \epsilon_j^{(2)} \in \C$.  It has mass $M_j$.

In general, the bending angle is expressed as an integral expressed linearly in terms of the mass distribution, see ~\cite[Equation 3.57]{PETTERS}. In particular, with point masses, the bending angle decomposes to a sum of bending angles, one for each point mass.  Each is computed as in the one mass system: 
\begin{align*}
	\tilde{\alpha}_j=\frac{4GM_j}{c^2 D_L \bar{\theta_j} } \qquad \mbox{where} \qquad \theta_j = \theta - \epsilon_j.
\end{align*}
We obtain
\begin{align*}\label{first}
	\beta=\theta-\sum_{j=1}^n \alpha_j = \theta-\frac{D_{LS}}{D_S D_L}\sum_{j=1}^n\frac{4GM_j}{c^2 \left(\bar{\theta}-\bar{\epsilon_j}\right)}.
\end{align*}

\noindent Letting 
\begin{align*}
w=\beta, \quad  z=\theta, \quad  z_j=\epsilon_{j}, \quad \mbox{and} \quad  \sigma_j=\frac{D_{LS}}{D_S D_L}\cdot\frac{4GM_j}{c^2}
\end{align*}
 gives
\begin{equation}\label{ORIGINAL LENS}
	w=z-\sum_{j=1}^n\frac{\sigma_j}{\bar{z}-\bar{z_j}}.
\end{equation}
Equation~\eqref{ORIGINAL LENS} requires the assumption that the center of mass is the origin, i.e. $\sum \sigma_j z_j = 0$. A translation by $w$ allows us to fix the position of the light source at the origin and vary the location of the center of mass. This simplifies Equation~\eqref{ORIGINAL LENS} to Equation~\eqref{LENS EQUATION}.

\subsection{Wilmshurst's Conjecture}
\index{Wilmshurst's Conjecture}
Let us finish our notes with an open problem that can be explored by undergraduates.
In \cite{WILMS}, Wilmshurst considered equations of the form
\begin{align}\label{EQN_WILMS}
p(z) = \bar{q(z)}
\end{align}
where $p(z)$ and $q(z)$ are polynomials of degree $n$ and $m$, respectively.
By conjugating the equation, if necessary, one may suppose $n \geq m$.  If $m=n$, then one
can have infinitely many solutions (e.g. $p(z) = z^n = q(z)$), but once $n > m$
Wilmshurst showed that there are finitely many solutions.  He conjectured that
the number of solutions to (\ref{EQN_WILMS}) is at most $3n-2 + m(m-1)$.

Unfortunately, this conjecture is false!  Counterexamples were  found when
$m=n-3$ by Lee, Lerario, and Lundberg  \cite{LLL}.  
They propose:

\begin{conjecture*}{\bf (Lee, Lerario, Lundberg)}\index{Lee}\index{Lerario}\index{Lundberg}
If $\deg(p(z)) = n$, $\deg(q(z))~=~m$, and $n > m$, then the number of solutions to $p(z) = \bar{q(z)}$
is bounded by  $2m(n-1) + n$.
\end{conjecture*}
\noindent
Note that this conjectured bound is not intended to be sharp.  For example,
Wilmshurst proved his conjecture in the case that $m=n-1$, providing
a stronger bound in that case  \cite{WILMS}.

This problem was further studied using certified numerics by Hauenstein,
Lerario, Lundberg, and Mehta \cite{HLLM}.\index{Hauenstein}\index{Mehta}  Their
work provides further evidence for this conjecture.

\begin{question*} Can techniques from complex dynamics be used to prove this conjecture?
\end{question*}

\bibliographystyle{plain}
\bibliography{complex_dynamics.bib}
\printindex

\end{document}